\documentclass[12pt]{amsart}
\usepackage{amsmath,amssymb}
\usepackage{color}
\usepackage[dvips]{graphicx}
\definecolor{gray}{gray}{.5}
\definecolor{gray2}{gray}{.9}
\def\gb{\makebox(0,0){\rotatebox{135}{\scalebox{1.5}{\color{gray}$\blacksquare$}}}}

\newcommand{\beqn}{\begin{eqnarray}}
\newcommand{\eeqn}{\end{eqnarray}}
\newcommand{\ep}{\varepsilon}

\newcommand{\C}{\mathbb{C}}

\newcommand{\Z}{\mathbb{Z}}

\newcommand{\sgn}{\mathrm{sgn}}
\newcommand{\cell}{\color{gray}\rule{3mm}{3mm}}
\newcommand{\Cell}{\color{gray}\rule{7mm}{7mm}}
\newcommand{\Cellb}{\color{gray2}\rule{8mm}{8mm}}
\newcommand{\Cellc}{\color{gray}\rule{4mm}{4mm}}
\newtheorem{thm}{Theorem}
\newtheorem{lem}{Lemma}
\newtheorem{prop}{Proposition}
\newtheorem{cor}{Corollary}
\newtheorem{Def}{Definition}
\newtheorem{Exam}{Example}

\begin{document}

\title{Excited Young diagrams and equivariant Schubert calculus}

\author{Takeshi Ikeda and Hiroshi Naruse}

\begin{abstract}
We describe the torus-equivariant cohomology ring of isotropic
Grassmannians by using a localization map to the torus
fixed points.
We present two types of formulas for equivariant 
Schubert classes 
of these homogeneous spaces.
The first formula involves combinatorial 
objects which we call ``excited Young diagrams'' and the 
second one is written in terms of factorial Schur
$Q$- or $P$-functions. As an application, we give
a Giambelli-type formula for the equivariant Schubert
classes.
We also give 
combinatorial and Pfaffian formulas for the multiplicity of a singular
point in a Schubert variety.
\end{abstract}

\maketitle

\section{Introduction}\label{Intro}
\setcounter{equation}{0}

In this paper, we give explicit descriptions 
of the Schubert classes
in the (torus) equivariant cohomology ring
of the Grassmannians as well as the maximal isotropic Grassmannians of both
symplectic and orthogonal types.
Our main results express the image of an equivariant Schubert class
under the {\it localization map\/} 
to the torus fixed points.

Now let us fix some notation.
Let $G$ be a complex semisimple connected  
algebraic group.
Choose a maximal torus $T$ of $G$ and a Borel subgroup $B$ containing $T.$
Let $P$ be a maximal parabolic subgroup
of $G$ containing $B.$ 
We are interested in the (integral) $T$-equivariant cohomology 
ring $H_T^*(G/P)$ of the homogeneous space $G/P.$ 
The equivariant Schubert classes are parametrized by
the set $W^P$ of minimal length representatives 
for $W/W_P,$ where $W$ is the Weyl groups of $G$
and $W_{P}$ is the parabolic subgroup associated
to $P.$
The set $W^P$ also parametrizes the $T$-fixed 
points $(G/P)^T$ in $G/P.$
In fact if we put $e_v=vP\;(v\in W^{P})$
then $(G/P)^{P}=\{e_{v}\}_{v\in W^{P}}.$
Let $B_{-}$ denote the opposite Borel subgroup 
such that $B_{-}\cap B=T.$
Define the Schubert variety $X_{w}$ 
associated with the element $w\in W^P,$  
to be the closure 
of $B_{-}$-orbit $B_{-}e_w$ of $e_{w}.$
Note that the codimension of $X_w$ in $G/P$ is 
$\ell(w)$, the length of $w,$ 
and $e_v\in X_w$ if and only if $w\leq v,$ where
$\leq$ is the partial order on $W^P$ induced by the Bruhat-Chevalley
ordering of $W.$ 
Since $X_w$ is a $T$-stable subvariety in $G/P,$ it induces
a $T$-equivariant fundamental class, the equivariant 
Schubert class, denoted by $[X_w]\in H_T^{2\ell(w)}(G/P).$
Our main goal is to describe $[X_w]$ explicitly.

In this paper, we consider $G/P$ in the following list:
\begin{itemize}
\item {\bf Type} $A_{n-1}$: $SL(n)/P_d\;(1\leq d\leq n),$
\item {\bf Type} $B_n$: $SO(2n+1)/P_n,$
\item {\bf Type} $C_n$: $Sp(2n)/P_n,$
\item {\bf Type} $D_n$: $SO(2n)/P_{d}\;(d=n-1,n)$
\end{itemize}
where we denote by 
$P_d$ the maximal parabolic subgroup
associated to the $d$-th simple root
(the simple roots being indexed as in \cite{Bo}).
It is well known that the space $SL(n)/P_d$ can be identified with 
the Grassmannian $\mathcal{G}_{d,n}$ of $d$-dimensional 
subspaces in $\C^n.$
Any other $G/P$ in the above list is 
a maximal isotropic Grassmannian
with respect to an orthogonal
or symplectic form  
(see Section \ref{BCD} for details).

Our description is based on
the ring homomorphism
$$
\iota^*: H^*_T(G/P)\longrightarrow
H^*_T((G/P)^T)=\bigoplus_{v\in W^P}H_T^*(e_v),
$$
induced by the inclusion $\iota: (G/P)^T\hookrightarrow G/P.$
This $\iota^*$ is known to be injective and called the localization map.
Each summand $H_T^*(e_v)$ is canonically isomorphic 
to the symmetric algebra $S=Sym_\Z(\hat{T})$ of the character group
$\hat{T}$ of the torus $T.$
Thus the equivariant Schubert class $[X_w]$ is 
described by a list $\{[X_w]|_{v}\}_{v\in W^{P}}$ of
polynomials in $S,$ where $[X_w]|_{v}$ denote
the image of 
the equivariant Schubert class $[X_w]$
under the homomorphism 
$\iota_v^*:H_{T}^{*}(G/P)\longrightarrow
H_{T}^{*}(e_{v})$ induced by 
the inclusion $\iota_v: \{e_v\}\hookrightarrow G/P.$

In type $A_{n-1}$ case (cf. Theorem \ref{AFac}), Knutson and Tao \cite{KT}
discovered that $[X_w]|_v$ can be identified
with a suitably specialized `factorial' Schur function,
a multi-parameter deformation of Schur function (see Section \ref{SectFac} for 
the definition). Their argument uses a remarkable 
vanishing property of the factorial Schur function 
(cf. Proposition \ref{vanishing}). 
By a totally different method,
Lakshmibai, Raghavan, and Sankaran \cite{LRS}
showed the same result, although they do not 
state it explicitly in terms of factorial Schur function. 
In fact, they start from 
a combinatorial expression for $[X_w]|_v$ 
in terms of a set of {\it non-intersecting paths\/},
which comes from a detailed analysis 
of Gr\"obner basis of the defining ideal
of the Schubert variety due to 
Kreiman and Lakshmibai \cite{KL},
and Kodiyalam and Raghavan \cite{KR},
and then rewrite the expression into a ratio of 
some determinants, which is a form of factorial Schur function.

Type $C_n,$ 
the case of Lagrangian Grassmannian,
was studied in a paper \cite{Ik} by the first named
author, where
$[X_w]|_v$ is expressed in terms 
of factorial Schur $Q$-function defined by Ivanov \cite{Iv}. 
The proof is a comparison 
of Pieri-Chevalley type recurrence relations
for both $[X_w]|_v$ and the factorial Schur $Q$-function.
This strategy of identification goes well for
other $G/P$ in our list above.
Actually, we prove in this paper the analogous
result for types $B_n$ and $D_n$, the orthogonal
Grassmannian, i.e., 
we present a formula for $[X_w]|_v$ in terms 
of factorial Schur $P$-function for 
these spaces (Theorem \ref{BCD-closed}).

As an application of these
formulas, we obtained an Giambelli-type formula (Corollary \ref{Giam})
for isotropic Grassmannians
that expresses an arbitrary equivariant Schubert class
as a Pfaffian of Schubert classes 
associated with the `two-row' (strict) partitions.
This formula is an equivariant analogue of 
the Giambelli formula due to Pragacz \cite{Pr}
in the case of ordinary cohomology.

Another type of formulas (Theorem \ref{EYD}, Theorem \ref{BCD}) 
we discuss in the present
paper involves combinatorial objects,
which we call {\it excited Young diagrams\/}, EYDs for short.
The idea of EYDs was inspired by the work \cite{LRS} of 
Lakshmibai, Raghavan, and Sankaran mentioned above.
These formulas have a `positive' nature
in the sense that it is expressed as 
a sum over a set of EYDs with each summand
being a products of some positive roots.

For an exposition,
here we consider $\mathcal{G}_{d,n}$.
It is well-known that the set $W^{P_d}$ is
parametrized by the set of
partition $\lambda=(\lambda_1,\ldots,\lambda_d)$
such that $n-d\geq \lambda_1\geq\cdots\geq \lambda_d\geq 0,$
or equivalently the Young diagrams contained in the rectangular
of shape $d\times (n-d).$
Let us denote by $D_\lambda$ the Young diagram of $\lambda.$
Suppose $w,v$ be elements of $W^{P_d}$ such that $e_v\in X_w.$
Let $\lambda,\mu$ be the corresponding partitions 
for $w,v$ respectively. Then
we have $D_\lambda\subset D_\mu.$
Our formula express $[X_w]|_{v}$
as a weighted sum over a set $\mathcal{E}_\mu(\lambda)$
(see Subsection \ref{DefEYD} for the definition) of 
subsets of $D_\mu.$
For example, let $\lambda=(3,2),\,\mu=(4,4,3,1).$
Some typical elements in $\mathcal{E}_\mu(\lambda)$ are illustrated below.
Here we depict the Young diagrams 
in Russian style.

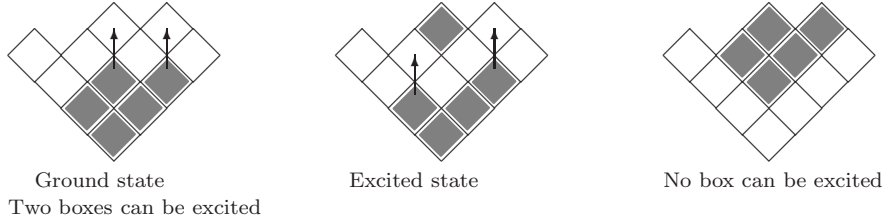
\begin{picture}(120,100)
\put(10,60){\line(1,-1){40}}
\put(20,70){\line(1,-1){40}}
\put(40,70){\line(1,-1){30}}
\put(50,80){\line(1,-1){30}}
\put(70,80){\line(1,-1){20}}
\put(10,60){\line(1,1){10}}
\put(20,50){\line(1,1){30}}
\put(30,40){\line(1,1){40}}
\put(40,30){\line(1,1){40}}
\put(50,20){\line(1,1){40}}
\put(50,30){\gb}
\put(40,40){\gb}
\put(50,50){\gb}
\put(60,40){\gb}
\put(70,50){\gb}
\put(50,55){\vector(0,1){15}}
\put(70,55){\vector(0,1){15}}
\put(20,10){\tiny{Ground state}}
\put(10,0){\tiny{Two boxes can be excited}}
\end{picture}
\begin{picture}(120,100)
\put(10,60){\line(1,-1){40}}
\put(20,70){\line(1,-1){40}}
\put(40,70){\line(1,-1){30}}
\put(50,80){\line(1,-1){30}}
\put(70,80){\line(1,-1){20}}
\put(10,60){\line(1,1){10}}
\put(20,50){\line(1,1){30}}
\put(30,40){\line(1,1){40}}
\put(40,30){\line(1,1){40}}
\put(50,20){\line(1,1){40}}
\put(50,30){\gb}
\put(40,40){\gb}
\put(50,70){\gb}
\put(60,40){\gb}
\put(70,50){\gb}
\put(40,45){\vector(0,1){15}}
\put(70,55){\vector(0,1){15}}
\put(15,10){\tiny{Excited state}}
\end{picture}
\begin{picture}(100,100)
\put(10,60){\line(1,-1){40}}
\put(20,70){\line(1,-1){40}}
\put(40,70){\line(1,-1){30}}
\put(50,80){\line(1,-1){30}}
\put(70,80){\line(1,-1){20}}
\put(10,60){\line(1,1){10}}
\put(20,50){\line(1,1){30}}
\put(30,40){\line(1,1){40}}
\put(40,30){\line(1,1){40}}
\put(50,20){\line(1,1){40}}
\put(50,50){\gb}
\put(40,60){\gb}
\put(50,70){\gb}
\put(60,60){\gb}
\put(70,70){\gb}
\put(10,10){\tiny{No box can be excited}}
\end{picture}

\noindent 
Let us imagine that 
each diagram labels a `quantum state'.
Each box can be {\it excited\/}
to the upper (Northern) space 
if the neighbors in North, North-East, and North-West directions
are all unoccupied (see Subsection \ref{DefEYD} for the precise definition
of {\it excitation\/}).
In this case, there are nine excited states
obtainable by applying successive excitations 
starting from the {\it ground state\/} i.e. the element $D_\lambda.$

It should be mentioned that
the notion of EYDs and its shifted analogue 
were introduced by Kreiman \cite{Kr1,Kr2}
independently to us.  
Kreiman \cite{Kr1} proved that
the set of non-intersecting paths appeared
in \cite{KL}, \cite{KR}, \cite{LRS} is naturally 
bijective to the set of EYDs.
He also presented a combinatorial formula of $[X_w]|_v$
for type $C_n$ case 
in terms of {\it shifted analogue\/} of EYDs
by using a result 
by Ghorpade and Raghavan \cite{GR} 
analogous to \cite{KL}, \cite{KR}.
In this paper, we present a different proof for these results
without using the Gr\"obner machinery mentioned above.
An advantage of our method is 
that we can apply the same argument 
to isotropic Grassmannians of 
orthogonal type, for which 
no explicit description of the Gr\"obner 
basis is known. 
In order to deal with the case of `even' orthogonal Grassmannian,
we introduce another variant of shifted EYDs.

Our method begins by 
identifying the localized classes $[X_w]|_v$
as $\xi$-functions defined by 
Kostant and Kumar (\cite{KK}).
Then we can make use of a well-known 
formula (Proposition \ref{AJS_B}, cf. \cite{Bi},\cite{AJS})  
that expresses arbitrary $\xi$-function as a sum over a set of `reduced subwords'.
And then use 
a theory by Stembridge \cite{St1}
on {\it fully commutative} elements in 
Coxeter groups.
The theory enable us to establish  
a natural bijection
between the set of reduced sub words appearing in 
the sum formula and a certain set of EYDs.
This bijection is a technical heart of 
our proof of Theorem \ref{EYD} and Theorem \ref{BCD}.

It is known that the multiplicity $m_v(X_w)$ at $e_v$ in $X_w$ is closely related to $[X_w]|_v.$
Such a multiplicity 
has been studied in detail by many authors (see \cite{BL}).
By our combinatorial formula, we can express
$m_v(X_w)$ as the number of elements 
in a certain set of excited Young diagrams.
And also, we can obtain a closed formula
for $m_v(X_w),$ which is a specialization
of a factorial Schur function.
This leads to a
Pfaffian formulas for the multiplicity of a singular
point in a Schubert variety.

In Section \ref{XiKK}, we explain the relation 
between the polynomial $[X_w]|_v$ and the $\xi$-function.
Some fundamental properties of $\xi$-functions are 
presented for the later use.
We discuss the case of $\mathcal{G}_{d,n}$ in
Sections \ref{SectEYD}, \ref{PfEYD}, \ref{SectFac}.
We present the combinatorial formula in Section \ref{SectEYD}
and the proof is given in Section \ref{PfEYD}.
We give the closed formula in Section \ref{SectFac},
which can be read independently 
from the preceding two sections. 
The isotropic Grassmannians 
are treated in a 
parallel manner in 
Sections \ref{SectBCD}, \ref{PfBCD}, \ref{SectQP}.
In Section \ref{SectMult} we discuss  
some application for $m_v(X_w).$
In Section \ref{SectLPM}
we discuss a relation 
between two types of formula
(Theorem \ref{BCD} and \ref{BCD-closed}) 
by using a Gessel-Viennot 
type argument.
\bigskip

{\bf Acknowledgments.} The authors thank
M. \:Ishikawa, S.\:Okada, and 
M. \:Shimozono for valuable discussions.

\section{$\xi$-functions of Kostant and Kumar}\label{XiKK}
\setcounter{equation}{0}
In this section we introduce the family of functions
$\xi^w$ for $w\in W$ defined by Kostant and Kumar.
By virtue of the result of Arabia \cite{Ar}, 
we can identify $\xi^w\,(w\in W^P)$ with the equivariant 
Schubert class $[X_w]$ in $H_T^*(G/P).$

Let $R^{+}$ denote the set of positive roots
with respect to $B.$ For $\beta\in R_{+},$
we denote by $\beta^\vee$ its dual coroot and
$s_\beta\in W$ the reflection
corresponding to $\beta.$ 
Let $\{\alpha_1,\ldots,\alpha_r\}$ be the set of simple 
roots in $R^{+}$.
The Weyl group $W$ is a Coxeter group generated by the 
set of simple reflections $\{s_1,\ldots,s_r\},$
where $s_i=s_{\alpha_i}.$
In particular, we can talk of
the length $\ell(w)$ of any element $w\in W.$
For a simple reflection $\alpha$ we denote by $\varpi_\alpha$ 
the corresponding fundamental weight.

\begin{prop}(Kostant and Kumar \cite{KK})\label{xi} There exist a family
of functions
$
\xi^w:W\longrightarrow S
$
for $w\in W$ 
with the following properties: 
\begin{enumerate}
\item $\xi^w(v)$ equals zero unless $w\leq v,$
\item $\xi^w(w)=\prod_{\alpha\in R_+\cap wR_{-}}\alpha,$
\item $\xi^{e}(v)=1$ for all $v\in W,$ where $e$ is the identity of $W,$ 
\label{identity}
\item If $\alpha$ is a simple root, then for all $v\in W,$
$$\xi^{s_\alpha}(v)=\varpi_\alpha-v(\varpi_\alpha),$$
\item If $\alpha$ is a simple root, then
\beqn
(\xi^{s_\alpha}-\xi^{s_\alpha}(w))\xi^{w}=
\sum_{w\overset{\beta}{\to} w'}\langle w(\varpi_\alpha),\beta^{\vee}\rangle\xi^{w'},
\label{PieriXi}
\eeqn
where we use the notation $w\overset{\beta}{\to} w'$ to indicate
$w'=s_\beta w $ for some $\beta\in R_{+}$ and $\ell(w')=\ell(w)+1,$
\item Each $\xi^w(v)$ with $v\in W$ is homogeneous of degree $\ell(w).$
\end{enumerate}
\end{prop}

{\it Remark.} We use notation for 
$\xi$-functions in Kumar's book,
which is different from the one
used in \cite{KK}.

\bigskip

The function $\xi^w$ is directly related to 
the object of our main interest. 

\begin{prop}\label{A} Let $w\in W^P.$ We have
$$
[X_w]|_{v}=\xi^w(v) 
$$
for $v\in W^P.$
\end{prop}
{\it Proof.}
Arabia \cite{Ar} proved the result
for the flag variety $G/B$
(see Graham's paper \cite{Gr} for more 
information).
For the parabolic case,
the reader can consult Kumar's book \cite{Ku2}.
\hfill $\square$

\bigskip

Let $P=P_d$ denote the maximal parabolic 
subgroup associated to $\alpha_d.$
Note that we have for $w\in W^{P}$ 
\beqn
(\xi^{s_d}-\xi^{s_{d}}(w))\xi^{w}=
\sum_{w'\in W^{P},\;w\overset{\beta}{\to} w'}\langle w(\varpi_\alpha),\beta^{\vee}\rangle\xi^{w'}.
\label{ParaXi}
\eeqn
The above relation involves only $\xi^w\;(w\in W^{P}).$
This reflects
the fact that the Schubert classes 
$[X_w]\;(w\in W^{P})$ form an
$S$-basis of $H_T^*(G/{P}),$
considered a sub $S$-algebra of $H_T^*(G/B)$
via the projection $G/B\rightarrow G/{P}.$ 
Note that if $w\in W^P$, $\xi^{w}$ 
is $W_P$-invariant in the sense that 
$$
\xi^w(vu)=\xi^w(v)\quad
\mbox{for all}\quad
u\in W_P.
$$

\bigskip

We can make use of the following formula:
\begin{prop}\label{AJS_B}(\cite{AJS}, \cite{Bi})
Let $w,v\in W$ such that $w\leq v.$
Fix a reduced expression $s_{i_1}\cdots s_{i_k}$ for $v$.
Put 
\beqn
\beta_t=s_{i_1}\cdots s_{i_{t-1}}(\alpha_{i_t})\quad
\mbox{for}\quad 1\leq t\leq k.
\label{beta}
\eeqn
Then
\beqn
\xi^w(v)=\sum_{j_1,\ldots,j_s}\beta_{j_1}\cdots\beta_{j_s},
\label{AJS,B}
\eeqn
where the sum is over
all sequences $1\leq j_1<\cdots<j_s\leq k$
such that $s_{i_{j_1}}\cdots s_{i_{j_s}}$ is 
a reduced expression for $w.$
\end{prop}

Although the above formula is explicit, 
it still requires a lot of
calculations to get a concrete expression for $\xi^w(v)$ in general. 
If $w$ is an element
of $W^P$ for classical $G$ and $P$ in our list (cf. Section \ref{Intro}),
then we can give a nice combinatorial   
interpretation for the right hand side of (\ref{AJS,B}).

\section{Excited Young diagrams}\label{SectEYD}
\setcounter{equation}{0}
We fix positive integers $n,d$ such that $1\leq d\leq n.$ 
In this section, we
give a combinatorial formula (Theorem \ref{EYD}) 
for the restriction of the
equivariant Schubert classes
in the Grassmannian $\mathcal{G}_{d,n}$
to any torus fixed points.

\subsection{Schubert variety of $\mathcal{G}_{d,n}$}
Let $w\in S_n$ be a {\it Grassmannian permutation\/}, i.e.,
$$
w(1)<\cdots<w(d),\quad
w(d+1)<\cdots<w(n).
$$
When we identify the space $SL(n)/P_d$ 
with the Grassmannian $\mathcal{G}_{d,n}$ of $d$-dimensional
subspaces of $\C^n,$
the Schubert variety associated with $w$ 
is given by 
$$
X_w=\{V\in \mathcal{G}_{d,n}\;|\; \dim(V\cap F_{n-w(d-i+1)+1})\geq i\quad\mbox{for}\quad
1\leq i\leq d\},
$$
where $F_i=\mathrm{span}_\C\{\pmb{e}_{n-i+1},\ldots,\pmb{e}_n\}$
is the $i$-plane spanned by 
the last $i$ vectors 
in the standard $T$-basis $\pmb{e}_1,\ldots,\pmb{e}_n$ of $\C^n.$
 
\subsection{Partitions and Young diagrams}

Let $\lambda=(\lambda_1\geq \cdots\geq \lambda_r\geq 0)$
be a partition. To every partition $\lambda$ one associates
its Young diagram $D_\lambda$ which is 
the set of square boxes with coordinate 
$(i,j)\in \Z^2$ such that
$1\leq j\leq \lambda_i:$
$$
D_\lambda=\{(i,j)\in \Z^2\;|\;
1\leq i\leq r,\; 1\leq j\leq \lambda_i\}.
$$
The boxes in $D_\lambda$ are arranged 
in a plane 
with matrix-style coordinates. 
For example,
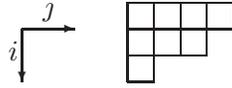
\begin{center}
\begin{picture}(85,45)
\put(50,10){\line(0,1){30}}
\put(60,10){\line(0,1){30}}
\put(70,20){\line(0,1){20}}
\put(80,20){\line(0,1){20}}
\put(90,30){\line(0,1){10}}
\put(50,10){\line(1,0){10}}
\put(50,20){\line(1,0){30}}
\put(50,30){\line(1,0){40}}
\put(50,40){\line(1,0){40}}
\put(10,30){\vector(1,0){20}}
\put(10,30){\vector(0,-1){20}}
\put(5,18){\small{$i$}}
\put(18,35){\small{$j$}}
\end{picture}
\end{center}
is the Young diagram $D_\lambda$ of $\lambda=(4,3,1).$

Let $\mathcal{P}_d$ denote the set of all partition
with length at most $d.$
If $n$ is an integer with $n\geq d,$
let $\mathcal{P}_{d,n}$ be the subset of $\mathcal{P}_d$
consisting of the elements 
whose largest part is less than or equal to $n-d.$
For any partition $\mu$, let
$\mathcal{P}_{\!\mu}$ denote the set of all partition
$\lambda$ such that $\lambda\leq \mu.$
In particular, if $\mu$ is the partition whose Young diagram is 
the $d\times (n-d)$ rectangle, then $\mathcal{P}_\mu$
is identical to $\mathcal{P}_{d,n}.$

\subsection{Grassmannian permutations and Young diagrams}
The set $W^{P_d}$ of minimal length coset representatives 
for $W/W_{P_d}$ with $W=S_n$ is 
identified with 
the set of all 
Grassmannian permutations.
Let $w\in W^{P_d}.$
We define a partition $\lambda=(\lambda_1,\ldots,\lambda_d)$ of 
$n$ by 
$$
\lambda_j=w(d-j+1)-d+j-1 \quad
(1\leq j\leq d).
$$
When considered as a Young diagram, $\lambda$ is contained 
in the rectangular shape 
$d\times (n-d).$
Note that we have 
$
\ell(w)=|\lambda|:=\sum_{j=1}^d\lambda_i,
$
where $\ell(w)$ is the length of $w.$

There is a convenient way to recover the Grassmannian permutation
from a Young diagram as follows.
Given a Young diagram $\lambda$ contained in 
the rectangle $d\times(n-d),$ write
a path along the boundary of the Young diagram
starting from the SW corner to the NE corner
of the rectangle.
We assign numbers to each 
arrow from $1$ to $n.$
For example, for $\lambda=(4,3,1,0)$ with $n=9,d=4$
we have the following picture:
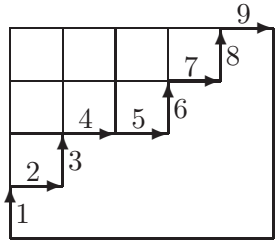
\setlength{\unitlength}{0.7mm}
\begin{center}
  \begin{picture}(60,50)
  \put(5,45){\line(1,0){40}}
  \put(5,35){\line(1,0){30}}
  \put(5,25){\line(1,0){10}}
  \put(5,15){\line(0,1){30}}
  \put(15,25){\line(0,1){20}}
  \put(25,25){\line(0,1){20}}
  \put(35,35){\line(0,1){10}}
  \put(5,5){\line(1,0){50}}
  \put(55,5){\line(0,1){40}}
  \thicklines
  \put(5,5){\vector(0,1){10}}
  \put(5,15){\vector(1,0){10}}
  \put(15,15){\vector(0,1){10}}
  \put(15,25){\vector(1,0){10}}
  \put(25,25){\vector(1,0){10}}
  \put(35,25){\vector(0,1){10}}
  \put(35,35){\vector(1,0){10}}
  \put(45,35){\vector(0,1){10}}
  \put(45,45){\vector(1,0){10}}
  %
  \put(6,8){\small{$1$}}
  \put(8,16){\small{$2$}}
  \put(16,18){\small{$3$}}
  \put(18,26){\small{$4$}}
  \put(28,26){\small{$5$}}
  \put(36,28){\small{$6$}}
  \put(38,36){\small{$7$}}
  \put(46,38){\small{$8$}}
  \put(48,46){\small{$9$}}
  %
  \end{picture}
\end{center}
If the assigned numbers of the vertical arrows are $i_1<\cdots<i_d$ and
those of the horizontal arrows are
$j_1<\cdots<j_{n-d},$ then the corresponding Grassmannian permutation
is 
$$
w=(i_1,\ldots,i_d,j_1,\ldots,j_{n-d}).
$$
Explicitly we have
$$
i_k=\lambda_{d-k+1}+k\quad
(1\leq k\leq d),\quad
j_k=-\lambda'_{k}+k+d\quad
(1\leq k\leq n-d),
$$
where $\lambda'$ is the conjugate of $\lambda.$
In the above example, we have $w=136824579.$

\subsection{Excited Young diagrams}\label{DefEYD}
Let $\lambda\leq \mu$ be partitions.
The Young diagram $D_{\lambda}$ of $\lambda$
is a subset of $D_{\mu}.$
Take an arbitrary subset $C$ of $D_{\mu}.$
Pick up a box $x\in C$ such that 
$x+(1,0),\,x+(0,1),\,x+(1,1)\in D_{\mu}\setminus C.$
Then set $C'=C\cup \{x+(1,1)\}\setminus \{x\}.$
\setlength{\unitlength}{0.4mm}
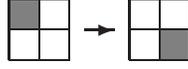
\begin{center}
\begin{picture}(65,25)
\put(5,15){\color{gray}\rule{4mm}{4mm}}
\put(55,5){\color{gray}\rule{4mm}{4mm}}
\multiput(5,5)(0,10){3}{\line(1,0){20}}
\multiput(5,5)(10,0){3}{\line(0,1){20}}
\multiput(45,5)(0,10){3}{\line(1,0){20}}
\multiput(45,5)(10,0){3}{\line(0,1){20}}
\thicklines
\put(30,15){\vector(1,0){10}}
\end{picture}
\end{center}
The procedure $C\rightarrow C'$ of changing $C$ into $C'$ is 
called an {\it elementary excitation\/}
occurring at $x.$
If a subset $S$ of $D_{\mu}$ is 
obtained from $C$ by applying elementary
excitations successively, i.e.,
there are sequence 
\beqn
C=C_{0}\rightarrow C_{1} \rightarrow \cdots \rightarrow C_{r-1}\rightarrow
C_{r}=S,\quad r\geq 0\label{seq}
\eeqn
of elementary excitations,
then we say that $S$ is an {\it excited state\/} of $C,$
or $S$ is obtained from $C$ by excitation.
Let $\mathcal{E}_{\mu}(\lambda)$ denote the set 
of all excited states of $D_{\lambda}.$

\begin{Exam}\label{EYDexample}{\rm
For example, let $\lambda=(3,1),\;\mu=(5,4,3).$
Then the set $\mathcal{E}_\mu(\lambda)$ 
consists of the following seven elements:}
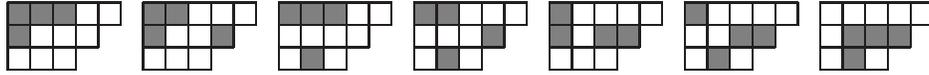
\setlength{\unitlength}{0.3mm}
\begin{center}
  \begin{picture}(400,40)
  \put(5,25){\cell}
  \put(15,25){\cell}
  \put(25,25){\cell}
  \put(5,15){\cell}
  \put(65,25){\cell}
  \put(75,25){\cell}
  \put(95,15){\cell}
  \put(65,15){\cell}
  \put(125,25){\cell}
  \put(135,25){\cell}
  \put(145,25){\cell}
  \put(135,5){\cell}
  \put(185,25){\cell}
  \put(195,25){\cell}
  \put(215,15){\cell}
  \put(195,5){\cell}
  \put(245,25){\cell}
  \put(265,15){\cell}
  \put(275,15){\cell}
  \put(245,15){\cell}
  \put(305,25){\cell}
  \put(325,15){\cell}
  \put(315,5){\cell}
  \put(335,15){\cell}
  \put(375,15){\cell}
  \put(385,15){\cell}
  \put(375,5){\cell}
  \put(395,15){\cell}
  \multiput(0,0)(60,0){7}{
  \put(5,35){\line(1,0){50}}
  \put(5,25){\line(1,0){50}}
  \put(5,15){\line(1,0){40}}
  \put(5,5){\line(1,0){30}}
  \put(5,5){\line(0,1){30}}
  \put(15,5){\line(0,1){30}}
  \put(25,5){\line(0,1){30}}
  \put(35,5){\line(0,1){30}}
  \put(45,15){\line(0,1){20}}
  \put(55,25){\line(0,1){10}}
  }%
  \end{picture} 
\end{center}
\end{Exam}

It is easy to see that the number $r$ in (\ref{seq}), the times of 
elementary excitations, is well-defined 
for $C\in \mathcal{E}_\mu(\lambda).$ 
In fact, if we define the {\it energy\/} $E(C)$ of $C\in \mathcal{E}_\mu(\lambda)$
by
$$E(C)=\sum_{(i,j)\in C}m(i,j)
-\sum_{(i,j)\in D_\lambda}m(i,j),\quad m(i,j)=\frac{1}{2}(i+j)
$$
then we have $E(C)=r.$

Let $\ep_1,\ldots,\ep_n$ 
be a standard basis of the 
lattice $L=\Z^n.$ 
The character group $\hat{T}$
is identified with a sub lattice of $L$ 
spanned by $\ep_i-\ep_{i+1}\;(1\leq i\leq n-1).$ 
Now we can state the first combinatorial formula.

\begin{thm}\label{EYD}(\cite{LRS},\cite{Kr1}) Let $w\leq v\in W^{P_d},$ 
and $\lambda\leq \mu\in \mathcal{P}_{d,n}$  the
corresponding partitions. Then we have
$$
[X_w]|_v=\sum_{C\in \mathcal{E}_\mu(\lambda)}
\prod_{(i,j)\in C}(\ep_{v(d+j)}-\ep_{v(d-i+1)}).
$$
\end{thm}

Proof of the theorem is given in the next section.

\begin{Exam}{\rm If $w=124735689,\,v=157923468\in W^{P_d}$ with $d=4,$ then 
$\lambda=(3,1),\mu=(5,4,3).$ We fill 
the boxes of $D_\mu$ with positive roots as follows:}
\setlength{\unitlength}{1.2mm}
\begin{center}
  \begin{picture}(60,50)
  \put(5,45){\line(1,0){50}}
  \put(5,35){\line(1,0){50}}
  \put(5,25){\line(1,0){40}}
  \put(5,15){\line(1,0){30}}
  \put(5,5){\line(0,1){40}}
  \put(15,15){\line(0,1){30}}
  \put(25,15){\line(0,1){30}}
  \put(35,15){\line(0,1){30}}
  \put(45,25){\line(0,1){20}}
  \put(55,35){\line(0,1){10}}
  \put(6,19){\small{$\ep_2\!-\!\ep_5$}}
  \put(6,29){\small{$\ep_2\!-\!\ep_7$}}
  \put(16,29){\small{$\ep_3\!-\!\ep_7$}}
  \put(16,19){\small{$\ep_3\!-\!\ep_5$}}
  \put(26,19){\small{$\ep_4\!-\!\ep_5$}}
  \put(26,29){\small{$\ep_4\!-\!\ep_7$}}
  \put(6,39){\small{$\ep_2\!-\!\ep_9$}}
  \put(16,39){\small{$\ep_3\!-\!\ep_9$}}
  \put(26,39){\small{$\ep_4\!-\!\ep_9$}}
  \put(36,39){\small{$\ep_6\!-\!\ep_9$}}
  \put(36,29){\small{$\ep_6\!-\!\ep_7$}}
  \put(46,39){\small{$\ep_8\!-\!\ep_9$}}
  \put(0,8){\small{$1$}}
  \put(9,46){\small{$2$}}
  \put(0,18){\small{$5$}}
  \put(19,46){\small{$3$}}
  \put(29,46){\small{$4$}}
  \put(0,28){\small{$7$}}
  \put(39,46){\small{$6$}}
  \put(0,38){\small{$9$}}
  \put(49,46){\small{$8$}}
  \end{picture}
\end{center}
{\rm Then our formula reads (cf. Example \ref{EYDexample}):}
\begin{small}
\beqn
[X_w]|_v&=&
(\ep_2-\ep_9)(\ep_3-\ep_9)(\ep_4-\ep_9)(\ep_2-\ep_7)
+(\ep_2-\ep_9)(\ep_3-\ep_9)(\ep_6-\ep_7)(\ep_2-\ep_7)\nonumber\\
&+&(\ep_2-\ep_9)(\ep_3-\ep_9)(\ep_4-\ep_9)(\ep_3-\ep_5)
+(\ep_2-\ep_9)(\ep_3-\ep_9)(\ep_6-\ep_7)(\ep_3-\ep_5)\nonumber\\
&+&(\ep_2-\ep_9)(\ep_4-\ep_7)(\ep_6-\ep_7)(\ep_2-\ep_7)
+(\ep_2-\ep_9)(\ep_4-\ep_7)(\ep_6-\ep_7)(\ep_3-\ep_5)\nonumber\\
&+&(\ep_3-\ep_7)(\ep_4-\ep_7)(\ep_6-\ep_7)(\ep_3-\ep_5).\nonumber
\eeqn
\end{small}
\end{Exam}

\section{Proof of Theorem \ref{EYD}}\label{PfEYD}
\setcounter{equation}{0}
\subsection{Fully commutative elements}
Let $W$ be a Coxeter group. 
An element $w$ in $W$ is {\it fully commutative\/} if 
any reduced expression for $w$ can be obtained from any
other by using only the Coxeter relations that
involve commuting generators. 
It is known that 
every element $w$ in $W^{P}$ for every $(G,P)$ in our list
(cf. Section \ref{Intro})  
is fully commutative (\cite{St1}, Theorem 6.1).

\subsection{Row reading expression for
a Grassmannian permutation}
Let $v$ be a Grassmannian permutation in $W^{P_d}$
and $\mu\in \mathcal{P}_{d,n}$ be the corresponding partition.
To each box $(i,j)\in D_\mu$ we fill in 
the simple reflection $s_{d-i+j}.$
For example, let $v=(3571246)\in S_{7}$ with $d=3.$
The corresponding 
partition is $\mu=(4,3,2)$ and we have
the following table:
\setlength{\unitlength}{0.5mm}
\begin{center}
  \begin{picture}(60,40)
  \put(5,35){\line(1,0){40}}
  \put(5,25){\line(1,0){40}}
  \put(5,15){\line(1,0){30}}
  \put(5,5){\line(1,0){20}}
  \put(5,5){\line(0,1){30}}
  \put(15,5){\line(0,1){30}}
  \put(25,5){\line(0,1){30}}
  \put(35,15){\line(0,1){20}}
  \put(45,25){\line(0,1){10}}
  \put(7,28){$s_3$}
  \put(17,28){$s_4$}
  \put(27,28){$s_5$}
  \put(37,28){$s_6$}
  \put(7,18){$s_2$}
  \put(17,18){$s_3$}
  \put(27,18){$s_4$}
  \put(7,8){$s_1$}
  \put(17,8){$s_2$}
  \end{picture} 
\end{center}
We read the entry of the boxes of the Young diagram $D_\mu$ 
from right to left starting from the bottom row to the top row
and form a word
\beqn
s_{i_1}\cdots s_{i_k} \quad (k=|\mu|),\label{rowread}
\eeqn
which gives a reduced expression for $v.$
For example we have
$$
v=3571246=s_2s_1\cdot s_4s_3s_2\cdot s_6s_5s_4s_3.
$$
We call the word given by (\ref{rowread}) 
the {\it row-reading word\/} of $v.$

The row reading procedure gives 
a bijective map
$$
\varphi: D_\mu\longrightarrow \{1,\ldots,k\},\quad
k=|\mu|.
$$
For example, if $\mu=(4,3,1)$ then
the map $\varphi$ is expressed by the following tableau:
\setlength{\unitlength}{0.5mm}
\begin{center}
  \begin{picture}(60,40)
  \put(5,35){\line(1,0){40}}
  \put(5,25){\line(1,0){40}}
  \put(5,15){\line(1,0){30}}
  \put(5,5){\line(1,0){20}}
  \put(5,5){\line(0,1){30}}
  \put(15,5){\line(0,1){30}}
  \put(25,5){\line(0,1){30}}
  \put(35,15){\line(0,1){20}}
  \put(45,25){\line(0,1){10}}
  \put(8,27){$9$}
  \put(18,27){$8$}
  \put(28,27){$7$}
  \put(38,27){$6$}
  \put(8,17){$5$}
  \put(18,17){$4$}
  \put(28,17){$3$}
  \put(8,7){$2$}
  \put(18,7){$1$}
  \end{picture}
\end{center}
Now let $C$ be an arbitrary subset of $D_\mu$ and 
$\varphi(C)=\{j_1,\ldots,j_r\}$ with $j_1<\cdots<j_r$ be the image
of $\{1,\ldots,k\}$ under the map $\varphi.$
Let $s_{i_j}$ be the simple reflection
assigned to the box $\varphi^{-1}(j)$ in $D_\mu$
for $1\leq j\leq k.$ 
Then we put
\beqn
w_C=s_{i_{j_1}}\cdots s_{i_{j_r}}.\label{w_C}
\eeqn
For example, if $C$ is the subset of $D_\mu$ 
indicated by the following gray boxes 
\setlength{\unitlength}{0.4mm}
\begin{center}
  \begin{picture}(60,40)
  \put(5,25){\color{gray}\rule{4mm}{4mm}}
  \put(5,15){\color{gray}\rule{4mm}{4mm}}
  \put(25,15){\color{gray}\rule{4mm}{4mm}}
  \put(5,35){\line(1,0){40}}
  \put(5,25){\line(1,0){40}}
  \put(5,15){\line(1,0){30}}
  \put(5,5){\line(1,0){20}}
  \put(5,5){\line(0,1){30}}
  \put(15,5){\line(0,1){30}}
  \put(25,5){\line(0,1){30}}
  \put(35,15){\line(0,1){20}}
  \put(45,25){\line(0,1){10}}
  \end{picture}
\end{center}
then $w_C=s_4s_2s_3.$
In particular, if $C=D_\mu$, then $w_C$ is nothing but 
the row-reading word of $v.$
Given another Grassmannian permutation
$w$ such that $w\leq v,$ define
$$
\mathcal{R}_v(w)=\{C\subset D_\mu\;|\; \sharp C=\ell(w),\;w_C=w\}.
$$
Note that if $\lambda$ is the partition corresponding to $w$,
then $D_\lambda\subset D_\mu$ and $D_\lambda\in \mathcal{R}_v(w).$

\begin{lem}\label{CC'}
Let $C,C'$ be subsets of a Young diagram $D_\mu$
such that $C'$ is obtained from $C$ by an elementary excitation from $C.$
If $C$ belongs to $\mathcal{R}_v(w)$
then we have $C'\in \mathcal{R}_v(w).$ 
\end{lem}
{\it Proof.}
We may assume that $C'$ is obtained from $C$ by an elementary
excitation
occurred in the $2\times 2$-square of the corner $(i,j).$ 
Set $k:=d+j-i.$
Consider the following regions in the diagram $D_\mu$:
$$
R:=\{(i,a)\;|\;j\leq a\leq \mu_i\}\cup
\{(i+1,b)\;|\;1\leq b\leq j+1\},
$$
$$
R_{\sharp}:= \{(i,a)\;|\;j+2\leq a\leq \mu_i\},\quad
R_{\flat}:= \{(i+1,b)\;|\;1\leq b\leq j-1\}.
$$
If we pick up the $i$-th and the $(i+1)$-th rows,
they look like
\begin{center}
\setlength{\unitlength}{0.7mm}
\begin{picture}(160,30)
\put(20,0){\line(1,0){60}}
\put(20,10){\line(1,0){100}}
\put(60,20){\line(1,0){60}}
\put(70,0){\line(0,1){20}}
\put(80,0){\line(0,1){20}}
\put(60,0){\line(0,1){20}}
\put(20,0){\line(0,1){10}}
\put(120,10){\line(0,1){10}}
\put(95,13){$R_{\sharp}$}
\put(35,3){$R_{\flat}$}
\put(0,3){$i+1$}
\put(62,23){$j$}
\put(70,23){$j+1$}
\put(115,23){$\mu_i$}
\put(5,13){$i$}
\put(63.5,13){$\bullet$}
\put(73.5,13){$\circ$}
\put(63.5,3){$\circ$}
\put(73.5,3){$\star$}
\end{picture}
\end{center}
It suffices to compare the sub-words
corresponding to $C\cap R$ and $C'\cap R$
(given by row reading procedure).
In $C$, we have $s_k$ at the position $(i,j)$ indicated above
by $\bullet.$
The positions indicated by $\circ$ are 
vacant by definition of the elementary excitation.
Then $C'$ is obtained by moving $s_k$ to the position indicated by 
$\star.$
Any simple reflection $s_l$ located in $C\cap R_\sharp$ (resp. $C\cap R_\flat$)
commutes with $s_k$
since $l\geq j+2$ (resp. $l\leq j-2$).
Hence the sub-word corresponding to the subset 
$C\cap R$ can be
rewritten into the sub-word corresponding to $C'\cap R$
using only the Coxeter relation
that involve commuting generators.
\hfill $\square$

\bigskip
 
\begin{cor} We have $\mathcal{E}_\mu(\lambda)\subset \mathcal{R}_v(w).$ 
\end{cor}
{\it Proof.} Use 
induction on energy $E(C)$ of $C\in \mathcal{E}_\mu(\lambda).$
The corollary is obvious from Lemma \ref{CC'}.
\hfill
$\square$

\bigskip

We would like to establish the following.
\begin{prop}\label{bij}
We have $\mathcal{E}_\mu(\lambda)=\mathcal{R}_v(w).$
\end{prop}

In order to prove Proposition \ref{bij},
it suffices to show the following.

\begin{lem}\label{RE}
Let $C\in \mathcal{R}_v(w)$ be such that $C\ne D_\lambda.$
Then there exists an element $C'\in \mathcal{R}_v(w)$
such that $C$ is obtained from $C'$ by a sequence of elementary
excitations.
\end{lem}
{\it Proof.}
Let $s_{i_{j_1}}\cdots s_{i_{j_r}}$
be the row-reading word for $w,$
and 
$s_{i_{k_1}}\cdots s_{i_{k_r}}$
be the word corresponding to $C.$
Let $a$ be such that $j_a\ne k_a$ and 
$
j_t=k_t
\; \mbox{for}
\; a< t\leq r.
$
We shall
compare the two words
$$
s_{i_{k_1}}\cdots s_{i_{k_a}},\quad
s_{i_{j_1}}\cdots s_{i_{j_a}}.
$$
Note that the element, say $w'$, expressed by the words
is fully commutative.
Put $t=i_{j_a}.$
Since 
$\{i_{k_1},\ldots,i_{k_a}\}$
is equal to $\{i_{j_1},\ldots,i_{j_a}\}$
as a multiset, 
we have $t\in \{i_{k_1},\ldots,i_{k_{a}}\}.$
Let $b$ be the largest 
index such that $i_{k_b}=t.$
Since $w'$ is fully commutative, the 
simple reflections adjacent to $s_t$, i.e. $s_{t+1},s_{t-1}$, 
can not appear in the subword
$$
s_{i_{k_{b+1}}}\cdots s_{i_{k_a}}.
$$ 
This implies a region $R$ indicated by the picture below
is unoccupied in the diagram $C:$
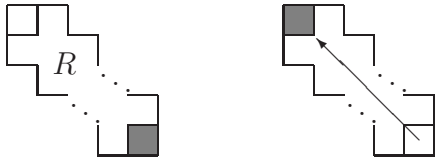
\begin{center}
\setlength{\unitlength}{0.4mm}
\begin{picture}(60,60)
\put(40,0){\color{gray}\rule{4mm}{4mm}}
\put(0,50){\line(1,0){20}}
\put(0,40){\line(1,0){10}}
\put(0,30){\line(1,0){10}}
\put(20,40){\line(1,0){10}}
\put(10,20){\line(1,0){10}}
\put(20,10){$\ddots$}
\put(30,20){$\ddots$}
\put(15,27){$R$}
\put(0,30){\line(0,1){20}}
\put(10,20){\line(0,1){10}}
\put(10,40){\line(0,1){10}}
\put(20,40){\line(0,1){10}}
\put(30,30){\line(0,1){10}}
\put(30,0){\line(0,1){10}}
\put(40,0){\line(0,1){10}}
\put(50,0){\line(0,1){20}}
\put(40,20){\line(1,0){10}}
\put(40,10){\line(1,0){10}}
\put(30,0){\line(1,0){20}}
\end{picture}
\hspace{1cm}
\begin{picture}(60,60)
\put(0,40){\color{gray}\rule{4mm}{4mm}}
\put(0,50){\line(1,0){20}}
\put(0,40){\line(1,0){10}}
\put(0,30){\line(1,0){10}}
\put(20,40){\line(1,0){10}}
\put(10,20){\line(1,0){10}}
\put(20,10){$\ddots$}
\put(30,20){$\ddots$}
\put(45,5){\vector(-1,1){34}}
\put(0,30){\line(0,1){20}}
\put(10,20){\line(0,1){10}}
\put(10,40){\line(0,1){10}}
\put(20,40){\line(0,1){10}}
\put(30,30){\line(0,1){10}}
\put(30,0){\line(0,1){10}}
\put(40,0){\line(0,1){10}}
\put(50,0){\line(0,1){20}}
\put(40,20){\line(1,0){10}}
\put(40,10){\line(1,0){10}}
\put(30,0){\line(1,0){20}}
\end{picture}
\end{center}
So we can move the box corresponding to 
$s_{k_b}$ as illustrated by the above 
picture to get a subset $C'$ of $D_\mu.$
Clearly $C'$ belongs to $\mathcal{R}_v(w).$ 
Since $C'$ has strictly smaller energy than $C,$ 
we see that $C'$
belongs to $\mathcal{E}_\mu(\lambda)$ by inductive hypothesis.
Now by the construction, $C$ is obtained from $C'$ by a sequence of 
elementary excitation. So by Lemma \ref{CC'}, implies 
$C$ is also a member of 
$\mathcal{E}_\mu(\lambda).$ 
\hfill $\square$

\subsection{$\beta$-sequences}

Fix a Grassmannian permutation $v,$
and let $\mu\in \mathcal{P}_{d,n}$ be the corresponding
partition. 
Let
$\varphi:D_\mu\longrightarrow\{1,2,\ldots,k\} $
be the row-reading map of $\mu.$
For $(i,j)\in D_\mu$ set $\beta_{i,j}=\beta_{\varphi(i,j)},$
where $\beta_t\:(1\leq t\leq k)$ are defined
in (\ref{beta}). 
\begin{lem}\label{betalem1}
For $(i,j)\in D_\mu$, we have 
$$
\beta_{i,j}=\ep_{v(d+j)}-\ep_{v(d-i+1)}.
$$
\end{lem}

{\it Proof.} To prove the lemma, we proceed by induction on $|\mu|=\ell(v).$
In case $\ell(v)=0,$ there is nothing to prove.
If $\ell(v)>0$, let 
$
v=s_{i_1}s_{i_2}\cdots s_{i_\ell}
$
be the row-reading expressions for $v.$
Setting $v'=s_{i_2}\cdots s_{i_\ell}$, which is the 
row-reading expression for $v$, and in particular we 
have $\ell(v')=\ell(v)-1.$
The corresponding shape $\mu'$ for $v'$ is 
obtained from $\mu$ by deleting a box
of position $(r,\mu_r)$, the left most one in the bottom row,
where $r$ is the number of rows in $\mu.$
Then we have
$$
\beta_{r,\mu_r}=\alpha_{i_1}=\ep_{d-r+\mu_r}-\ep_{d-r+\mu_r+1}.
$$
By definition of $\mu_r$, we have 
$\mu_r=v(d-r+1)-d+r-1$, so
$$
\beta_{r,\mu_r}=\ep_{v(d-r+1)-1}-\ep_{v(d-r+1)}.
$$
Then by Lemma \ref{corner} below, we have 
$
\beta_{r,\mu_r}=\ep_{v(d+\mu_r)}-\ep_{v(d-r+1)}.
$

Let $(i,j)\in \mu'.$ 
By the hypothesis of induction, we have
$$
\beta_{i,j}'=\ep_{v'(d+j)}-\ep_{v'(d-r+i)}.
$$
By definition of $\beta_{i,j}$, it is
easy to see that $\beta_{i,j}=s_{i_1}(\beta_{i,j}').$
Hence we have
$$
\beta_{i,j}=s_{i_1}(\beta_{i,j}')=s_{i_1}(\ep_{v'(d+j)}-\ep_{v'(d-r+i)})
=\ep_{s_{i_1}v'(d+j)}-\ep_{s_{i_1}v'(d-r+i)}.
$$
Since $s_{i_1}v'=v$, we have the lemma.
\hfill $\square$

\begin{lem}\label{corner} 
Let $r$ be the number of rows in $\mu.$
Then we have 
$$v(d-r+1)=v(d+\mu_r)+1.
$$
\end{lem} 
{\it Proof.} 
For a sequence $a_1<a_2<\cdots <a_m$ of integers,
we say $a_i$ is a ``gap'' if $a_i-a_{i-1}\geq 2.$ Note that
$\mu_r$ is the largest number such that  
$v(d+1)<\cdots<v(d+\mu_r)$ has no gap.
Since $v$ is a Grassmannian permutation,
the number $v(d+\mu_r)+1$ occurs 
as the smallest gap in $0<v(1)<\cdots<v(d).$
On the other hand, we have $\mu_d=\cdots=\mu_{r+1}=0$
and $\mu_r>0,$ which implies that the smallest gap
in $0<v(1)<\cdots<v(d)$ is $v(d-r+1).$
\hfill $\square$

\bigskip

{\it Proof of Theorem} \ref{EYD}. 
By Proposition \ref{AJS_B} together with Lemma \ref{betalem1}
we have
$$
[X_w]|_v=\sum_{C\in \mathcal{R}_v(w)}\prod_{(i,j)\in C}
(\ep_{v(d+j)}-\ep_{v(d-i+1)}).
$$
Then the theorem is immediate from Proposition \ref{bij}.
\hfill $\square$

\section{Factorial Schur functions}\label{SectFac}
\setcounter{equation}{0}
Our main goal in this section is to 
express $[X_w]|_v$ for $\mathcal{G}_{d,n}$ as a specialization 
of a factorial Schur function.
First we recall the definition of the factorial Schur
functions. 
\subsection{Definition of factorial Schur functions}
Let $x=(x_1,\ldots,x_d)$ be a finite sequence of variables
and let $a=(a_i)_{i=1}^\infty$, be an infinite sequence of parameters.
The {\it factorial Schur function} for a partition
$\lambda$ of length at most $d$ can be defined as follows.
Let 
$$
(z|a)^k=(z-a_1)(z-a_2)\cdots(z-a_k)
$$
for any $k\geq 0.$
Then we put
$$
s_{\lambda}^{(d)}(x|a)
=\frac{\det((x_j|a)^{\lambda_i+d-i})_{1\leq i,j\leq d }}
{\prod_{1\leq i<j\leq d}(x_i-x_j)}.
$$
This function is
actually a polynomial in $x_1,\ldots,x_n$
and $a_1,a_2,\ldots,a_{\lambda_1+d-1},$
homogeneous of degree $|\lambda|.$
In particular we have
$$
s^{(d)}_1(x|a)=x_1+\cdots+x_d-a_1-\cdots-a_d.
$$

For a partition $\lambda\in \mathcal{P}_d$, we define 
a $d$-tuple
$$
a_\lambda=(a_{\lambda_d+1},a_{\lambda_{d-1}+2},\ldots,a_{\lambda_1+d}).
$$

\begin{prop}\label{vanishing}(Vanishing property, cf. \cite{MS})
We have $s^{(d)}_\lambda(a_\mu|a)=0$ unless
$\mu\geq \lambda.$
\end{prop}

Let $\lambda\in \mathcal{P}_d,$ and take
sufficiently large $n$ such that $\lambda$ is contained in 
the $d\times (n-d)$ rectangle. 
Define $n$-tuple $(w(1),\ldots,w(n))$
by
$$
w(i)=\begin{cases}
\lambda_{d-i+1}+i & \mbox{for}\;1\leq i\leq d\\
-\lambda'_{i}+i   & \mbox{for}\; d<i\leq n
\end{cases},
$$
then $w=(w(1),\ldots,w(n))$ is a permutation
of $\{1,\ldots,n\}.$
If we set $w(i)=i$ for all $i>n$, then 
the infinite sequence $(w(1),w(2),\ldots)$ does not
depend on the choice of $n.$

\begin{lem} (cf. \cite{MS}) Let $\lambda\in \mathcal{P}_d.$ We have
$$
s_\lambda^{(d)}(a_\lambda|a)
=\prod_{(i,j)\in D_\lambda}(a_{w(d-i+1)}-a_{w(d+j)}).
$$
\end{lem}

We also need the following
Pieri-type formula:
\begin{lem}\label{PieriA} (cf. \cite{MS}) Let $\lambda\in \mathcal{P}_d.$ We have
$$
\left(s^{(d)}_{(1)}(x|a)
-s^{(d)}_{(1)}(a_\lambda|a)\right)
s^{(d)}_\lambda(x|a)
=\sum_{\nu}s_\nu^{(d)}(x|a),
$$
where the summation in the right hand side runs through
$\nu\in \mathcal{P}_d$ such that $\nu\geq \lambda$ and 
$|\nu|=|\lambda|+1.$
\end{lem}

\subsection{Closed formulas}
Let $\lambda\in \mathcal{P}_{n,d}.$
Recall that $s^{(d)}_\lambda(x|a)$ is a polynomial
in $x_1,\ldots,x_d$ and $a_1,\ldots,a_{n-1}.$
A specialization given by 
$a_i=\ep_i\,(1\leq i\leq n-1)$
is important for our geometric application below.
For $v\in W^{P_d},$ we define an $n$-tuple by
$$
x_v=(\ep_{v(1)},\ldots,\ep_{v(d)}).
$$

\begin{lem}\label{divA} Let $v\in W^{P_d}.$ We have
$$
[X_{s_d}]|_v=-
s_1^{(d)}(x_{v}|\ep_1,\ldots,\ep_{n-1}).
$$
\end{lem}
{\it Proof.} The $d$-th fundamental weight of $SL(n)$ is
given by
$\varpi_d=\ep_1+\cdots+\ep_d.$
By Proposition \ref{xi}, (4) and formula (\ref{P1}), 
the lemma follows.
\hfill $\square$

\begin{thm}\label{AFac}
(Knutson and Tao \cite{KT}, Lakshmibai, Raghavan, and Sankaran \cite{LRS}) 
Let $w\leq v\in W^{P_d},$
and $\lambda\leq \mu\in \mathcal{P}_{d,n}$ be the 
associated partitions.
Then we have
\beqn
[X_w]|_{v}
=(-1)^{\ell(w)}s_\lambda^{(d)}(x_{v}|\ep_1,\ldots,\ep_{n-1}).\label{AFacEq}
\eeqn
\end{thm}

{\it Proof.} Let $\mathcal{P}_\mu$ denote the set of all partition
$\lambda$ such that $\lambda\leq \mu.$
Consider the following system of equations for
the functions 
$F_\lambda=F_\lambda(\ep_1,\ldots,\ep_n)\;(\lambda\in \mathcal{P}_{\mu}):$
\beqn
\sum_{i=1}^d (\ep_{\lambda_{d-i+1}+i}-\ep_{\mu_{d-i+1}+i})
\cdot F_\lambda
=\sum_{\nu}F_{\nu}\quad
(\lambda\in \mathcal{P}_\mu\setminus \{\mu\}),\label{rec}
\eeqn
where the summation in the right hand side runs through
$\nu\in \mathcal{P}_{\mu}$ such that $\nu\geq \lambda$ and 
$|\nu|=|\lambda|+1.$ This
equation together with the
initial condition $F_\phi(\ep_1,\ldots,\ep_n)=1$
determine the functions 
$F_\lambda(\ep_1,\ldots,\ep_n)\;(\lambda\in \mathcal{P}_{\mu})$
uniquely.

Now we can see equation (\ref{ParaXi})
for $[X_w]|_v$ is identical to (\ref{rec}) above,
in view of the bijection $W^{P_d}\cong\mathcal{P}_{d,n}.$ 
On the other hand, the right hand side of (\ref{AFacEq})
satisfies 
the same equation from Proposition \ref{vanishing}
and Lemmas \ref{PieriA} and \ref{divA}.
In addition, both hand sides of (\ref{AFacEq})
satisfy the initial condition.
Hence the theorem follows.
\hfill $\square$

\section{Lagrangian and Orthogonal Grassmannians}\label{SectBCD}
\setcounter{equation}{0}
Fix a positive integer $n.$
We give here analogous theorems 
for the some other classical homogeneous
spaces $G/P$ in the following types:
\begin{itemize}
\item {\bf Type} $B_n$: $G=SO(2n+1,\C)$, $P=P_n,$
\item {\bf Type} $C_n$: $G=Sp(2n,\C)$, $P=P_n,$
\item {\bf Type} $D_{n}$: $G=SO(2n,\C)$, $P=P_{n-1},\;P_n,$
\end{itemize}
where we shall denote by 
$P_d$ the maximal parabolic subgroup
associated to the simple root $\alpha_d$
(the simple roots being indexed as in \cite{Bo}). 

In type $B_n$ (resp. type $C_n$), the variety $G/P$ can be 
identified with a closed subvariety of 
$\mathcal{G}_{n,2n+1}$ (resp. $\mathcal{G}_{n,2n}$) 
parametrizing the isotropic $n$-spaces
in $\C^{2n+1}$ (resp. $\C^{2n}$)
equipped with a non-degenerate symmetric (resp. skew symmetric) form.
Our space $G/P$ in type $C_n$ is also called
the Lagrangian Grassmannian.
For the even dimensional space $\C^{2n}$
equipped with a non-degenerate symmetric form,
the isotropic $n$-subspaces constitutes a 
union of two closed subvariety of $\mathcal{G}_{n,2n}$ 
each of which is isomorphic to 
$G/P,$ $P$ being one of $P_{n-1}$ and $P_{n}.$ 
Note that the
varieties $SO(2n+1,\C)/P_n,\,SO(2n+2,\C)/P_{d}\;(d=n,n+1)$
are all isomorphic (cf. \cite{Pr}),
although they are different
as $T$-spaces.

The goal of this section is to present a combinatorial
formulas for the restriction of equivariant 
Schubert classes to any torus fixed points
for the classical Grassmannians $G/P.$

\subsection{The set $W^P$}
First we fix some notation on a set $W^P$ which
parametrizes the Schubert classes, and 
the torus fixed points.

The character group
$\hat{T}$ of our torus is
a lattice with a standard basis 
$\{\ep_i\}_{i=1}^r,$ where $r=\mathrm{dim}(T).$

{\bf Type} $C_n$ ($G=Sp({2n},\C),\,P=P_n$):
We identify $W$ with a subgroup of $S_{2n}$
acts on the functionals
$\pm\ep_1,\ldots,\pm\ep_n.$
Denote $\ep_i,-\ep_i$ by $i$ and $\bar{i}$ respectively
and define a partial order
on the set 
$I_{n}=\{1,\ldots,n,\overline{n},\ldots,\overline{1}\}$
by 
$$
1<2<\cdots<n<\overline{n}<\cdots<\overline{2}<\overline{1}.
$$
Then we have 
$
W=\{w\in S_{2n}\;|\; w(\,\overline{i}\,)=\overline{w(i)} \},
$
where $a\mapsto \overline{a}$
is the involution on $I_{C_n}$
given by $\pm \ep_i\mapsto \mp\ep_i.$
Then  
$$
W^P=\{w\in W\;|\; w(1)<\cdots<w(n) \}.
$$
The simple roots are 
$$
\alpha_i=\ep_i-\ep_{i+1}\quad(1\leq i\leq n),
\quad
\alpha_n=2\ep_n,
$$
and corresponding the simple reflections are
$$
s_i=(i,i+1)(\overline{i+1},\overline{i})\quad
(1\leq i\leq n-1),\quad
s_n=(n,\overline{n}).
$$
We have
$$
\varpi_n=\ep_1+\cdots+\ep_n.
$$

{\bf Type} $B_n$ ($G=SO({2n+1},\C),\,P=P_n$):
Since $W$ is identical to the case $C_n$
as a Coxeter group,
the description of $W^P$ is the same
as Type $C_n.$
The simple roots are given by 
$$
\alpha_i=\ep_i-\ep_{i+1}\quad(1\leq i\leq n),
\quad
\alpha_n=\ep_n,
$$
and we have 
$$
\varpi_n=\textstyle{\frac{1}{2}}(\ep_1+\cdots+\ep_n).
$$

{\bf Type} $D_{n+1}$ ($G=SO({2n+2},\C),\,P=P_{n+1}$):
In this case, we have
$$
W=\left\{
w\in S_{2n+2}\;|
 w(\,\overline{i}\,)=\overline{w(i)},\quad
\sharp N(w)\;\mbox{is even}
\right\},
$$
where we set 
$N(w)=\{\,i\;|1\leq i\leq n+1,\;w(i)>n+1\}.$
We have
$$
W^P=\{w\in W\;|\; w(1)<\cdots<w(n+1) \}.
$$
The simple reflections are given by 
$$s_i=(i,i+1)(\overline{i+1},\overline{i})\quad
(1\leq i\leq n),\quad
s_{n+1}=(n,\overline{n+1})(n+1,\overline{n}),$$
and we have
$$
\alpha_i=\ep_i-\ep_{i+1}\quad(1\leq i\leq n),
\quad
\alpha_{n+1}=\ep_n+\ep_{n+1},
$$
and
\beqn
\varpi_{n+1}=\textstyle{\frac{1}{2}}(\ep_1+\cdots+\ep_{n+1}).\label{pi}
\eeqn

{\it Remark 1.} For type $D_{n+1}$, the result
for $P=P_{n}$ is obtained by simply 
replacing $\ep_{n+1}$ by $-\ep_{n+1}$.
So we only consider $P_{n+1}.$

\subsection{Strict partitions and shifted Young diagrams}
Next we need a description
of the set $W^P$ in terms of 
strict partitions.

Let $\lambda=(\lambda_1>\cdots>\lambda_r>0)$ be 
a strict partition.
Then the {\it shifted Young diagram} of $\lambda$
is the array of boxes into $r$-rows with $\lambda_i$ boxes
in the $i$-th row, such that each row is shifted
by one position to the right relative to 
the preceding row.
More explicitly, given a strict partition  
$\lambda_1>\lambda_2>\cdots>\lambda_r>0,$ the shifted
diagram of $\lambda$ is defined to be
$$
D'_\lambda:=\{(i,j)\in \Z^2\;|\;1\leq i\leq l,\;
i\leq j<\lambda_i+i\}.
$$
For example
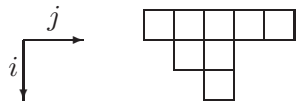
\begin{center}
\begin{picture}(100,55)
\put(50,40){\line(0,1){10}}
\put(60,30){\line(0,1){20}}
\put(70,20){\line(0,1){30}}
\put(80,20){\line(0,1){30}}
\put(90,40){\line(0,1){10}}
\put(100,40){\line(0,1){10}}
\put(50,50){\line(1,0){50}}
\put(70,20){\line(1,0){10}}
\put(60,30){\line(1,0){20}}
\put(50,40){\line(1,0){50}}
\put(10,40){\vector(1,0){20}}
\put(10,40){\vector(0,-1){20}}
\put(5,28){\small{$i$}}
\put(18,45){\small{$j$}}
\end{picture}
\end{center}
is the shifted Young diagram of $\lambda=(5,2,1).$

Define a partial order 
$\lambda\leq \mu$
if and only if $\lambda_i\leq \mu_i$ for all $i.$
Denote 
the particular element $\rho_n=(n,\ldots,2,1).$
Let $\mathcal{SP}_{\!\mu}$ denote the set of 
all strict partitions $\lambda$ such that $\lambda\leq \mu.$
The cardinality of the set $\mathcal{SP}_{\!\rho_n}$ is $2^n.$

\begin{prop}
Let $G/P$ as above. 
There is a natural
order-preserving bijection 
$$W^P\cong \mathcal{SP}_{\!\rho_n},$$
where if $\lambda\in \mathcal{SP}_{\!\rho_n}$
corresponds to $w\in W^P$ we have $\ell(w)=|\lambda|.$
\end{prop}

Explicitly, the bijections are 
given as follows:

{\bf Types} $B_n,C_n$: Let $w\in W^P.$ 
Since $w$ 
is a Grassmannian permutation of 
$2n$ letters (with $d=n$) 
we have the associated Young diagram, say $\Lambda.$ 
It is {\it symmetric\/} and
contained in an $n\times n$ square.
Then the associated element $\lambda\in \mathcal{SP}_{\!\rho_n}$
is the upper part
of the symmetric Young diagram including the diagonal.
Thus if $\Lambda=(\Lambda_1,\ldots,\Lambda_{n}),$ then 
$$
\lambda_i=\mathrm{max}\{\Lambda_i-i+1,0\}\quad (1\leq i\leq n).
$$

{\bf Type} $D_{n+1}$: Let $w\in W^P.$ 
Then the associated Young diagram, say $\Lambda$, 
is {\it symmetric\/} and
contained in an $(n+1)\times (n+1)$ square.
Then the associated element $\lambda\in \mathcal{SP}_{\!\rho_n}$
is the ``strictly'' upper part
of the symmetric Young diagram.
Thus if $\Lambda=(\Lambda_1,\ldots,\Lambda_{n+1}),$ then 
$$
\lambda_i=\mathrm{max}\{\Lambda_i-i,0\}\quad (1\leq i\leq n).
$$

\begin{Exam} 
{\rm Take $w_{1}=24\bar{3}\bar{1}$ for $D_{4}$ and 
$w_{2}=2\bar{3}\bar{1}$ for $B_{3}$ (or $C_{3}$).
Then the associated strict partitions are the same and given 
by $\lambda=(3,1).$}
\end{Exam}
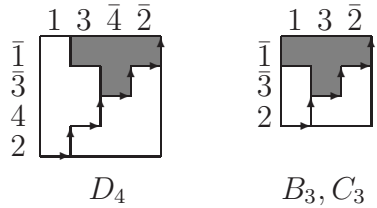
\setlength{\unitlength}{0.8mm}
\begin{center}
  \begin{picture}(65,40)
  \put(10,25){$\Cellc$}
  \put(15,25){$\Cellc$}
  \put(20,25){$\Cellc$}
  \put(15,20){$\Cellc$}
  \multiput(5,10)(5,5){4}{\vector(1,0){5}}
  \multiput(10,10)(5,5){4}{\vector(0,1){5}}
  \put(10,10){\line(1,0){15}}
  \put(25,10){\line(0,1){15}}
  \put(5,10){\line(0,1){20}}
  \put(5,30){\line(1,0){20}}
  \put(10,25){\line(0,1){5}}
  \put(10,25){\line(1,0){5}}
  \put(15,20){\line(0,1){5}}
  \put(0,10){$2$}
  \put(0,15){$4$}
  \put(0,20){$\bar{3}$}
  \put(0,25){$\bar{1}$}
  \put(6,31){$1$}
  \put(11,31){$3$}
  \put(16,31){$\bar{4}$}
  \put(21,31){$\bar{2}$}
  \put(45,25){$\Cellc$}
  \put(50,25){$\Cellc$}
  \put(55,25){$\Cellc$}
  \put(50,20){$\Cellc$}
  \multiput(45,15)(5,5){3}{\vector(1,0){5}}
  \multiput(50,15)(5,5){3}{\vector(0,1){5}}
  \put(50,15){\line(1,0){10}}
  \put(60,15){\line(0,1){10}}
  \put(45,15){\line(0,1){15}}
  \put(45,30){\line(1,0){15}}
  \put(45,25){\line(1,0){5}}
  \put(50,20){\line(0,1){5}}
  \put(41,15){$2$}
  \put(41,20){$\bar{3}$}
  \put(41,25){$\bar{1}$}
  \put(46,31){$1$}
  \put(51,31){$3$}
  \put(56,31){$\bar{2}$}
  \put(13,3){$D_{4}$}
  \put(45,3){$B_{3},C_{3}$}
\end{picture}
\end{center}
Note that $w_{1}=s_{3}s_{1}s_{2}s_{4}$ and $w_{2}=s_{3}s_{1}s_{2}s_{3}$
and these are the row-reading reduced expressions introduced in Section \ref{PfBCD}.

\subsection{Excited Young diagrams for shifted cases}
The notion of excitation of shifted Young diagrams 
is defined in the same way as in the case of ordinary 
Young diagrams if we 
introduce the idea of elementary excitation.

Let $\mu$ be a strict partition.
Given a subset $C$ of $D_{\mu}'.$
If a box $x\in C$ satisfies either of the following
conditions: 
\begin{enumerate}
\item $x=(i,j),\;i<j$ and $x+(1,0),\;x+(0,1),\;x+(1,1)\in D_{\mu}'\setminus C,$
\item $x=(i,i)$ and $x+(0,1),\;x+(1,1)\in D_{\mu}'\setminus C$
\end{enumerate}
then we set $C'=C\cup \{x+(1,1)\}\setminus \{x\}$
and call this procedure $C\rightarrow C'$ an elementary
excitation of type $I$ occurring at $x\in C.$
In general if a subset $S$ of $D_{\mu}'$ is obtained 
from $C$ by applying elementary excitations of type $I$ successively,
then we say that $S$ is an excited state of $C.$
Suppose we are given a strict partition $\lambda$ such that $\lambda\leq \mu.$
Let us denote by $\mathcal{E}_{\mu}^{I}(\lambda)$ 
the set of all excited states of $D_{\lambda}'.$

We also define the set $\mathcal{E}_{\mu}^{I\!I}(\lambda)$ 
consisting of the elements obtained 
from $C$ by a successive application of elementary excitations of ``type $I\!I$''
defined below.
Let $C$ be a subset of $D_{\mu}'$ and 
suppose we take a box $x\in C$ satisfying either of the following
conditions: 
\begin{enumerate}
\item $x=(i,j),\;i<j$ and $x+(1,0),\;x+(0,1),\;x+(1,1)\in D_{\mu}'\setminus C,$
\item $x=(i,i)$ and $x+(0,1),\;x+(1,1),\;x+(1,2),\;x+(2,2)\in D_{\mu}'\setminus C.$
\end{enumerate}
In the case of (1) we set $C'=C\cup \{x+(1,1)\}\setminus \{x\},$
and in the case of (2) $C'=C\cup \{x+(2,2)\}\setminus \{x\}.$
We call the procedure $C\rightarrow C'$  an elementary
excitation of type $I\!I$ occurring at $x\in C.$
Clearly we have $\mathcal{E}_{\mu}^{I\!I}(\lambda)\subset \mathcal{E}_{\mu}^{I}(\lambda).$

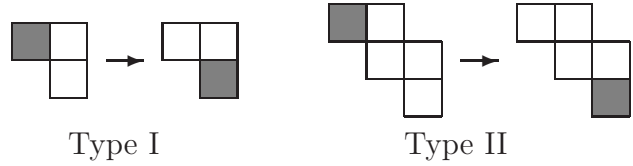
\setlength{\unitlength}{0.5mm}
\begin{center}
\begin{picture}(65,25)
\put(5,25){\color{gray}\rule{5mm}{5mm}}
\put(55,15){\color{gray}\rule{5mm}{5mm}}
\multiput(5,25)(0,10){2}{\line(1,0){20}}
\multiput(15,15)(10,0){2}{\line(0,1){20}}
\put(15,15){\line(1,0){10}}
\put(5,25){\line(0,1){10}}
\multiput(45,25)(0,10){2}{\line(1,0){20}}
\multiput(55,15)(10,0){2}{\line(0,1){20}}
\put(55,15){\line(1,0){10}}
\put(45,25){\line(0,1){10}}
\thicklines
\put(30,25){\vector(1,0){10}}
\put(20,0){Type I}
\end{picture}
\qquad
\begin{picture}(65,45)
\put(5,30){\color{gray}\rule{5mm}{5mm}}
\put(75,10){\color{gray}\rule{5mm}{5mm}}
\put(5,40){\line(1,0){20}}
\put(5,30){\line(1,0){30}}
\put(15,20){\line(1,0){20}}
\put(25,10){\line(1,0){10}}
\put(35,10){\line(0,1){20}}
\put(25,10){\line(0,1){30}}
\put(15,20){\line(0,1){20}}
\put(5,30){\line(0,1){10}}
\put(55,40){\line(1,0){20}}
\put(55,30){\line(1,0){30}}
\put(65,20){\line(1,0){20}}
\put(75,10){\line(1,0){10}}
\put(85,10){\line(0,1){20}}
\put(75,10){\line(0,1){30}}
\put(65,20){\line(0,1){20}}
\put(55,30){\line(0,1){10}}
\thicklines
\put(40,25){\vector(1,0){10}}
\put(25,0){Type II}
\end{picture}
\end{center}

\begin{Exam}\label{strictEYD}
{\rm Let $\lambda=(3,1),\;\mu=\rho_4.$ 
The set $\mathcal{E}_{\mu}^{I}(\lambda)$
consists of the following ten elements:}
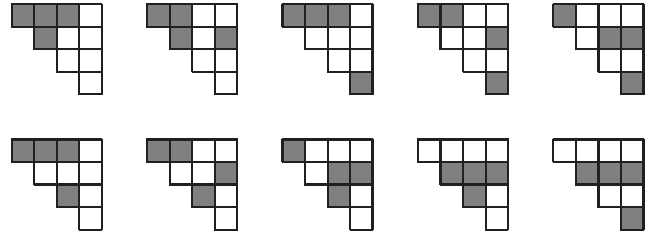
\setlength{\unitlength}{0.3mm}
\begin{center}
  \begin{picture}(300,120)
  \multiput(5,95)(60,0){5}{\put(0,0){\cell}}
  \multiput(25,95)(120,0){2}{\put(0,0){\cell}}
  \multiput(15,95)(60,0){4}{\put(0,0){\cell}}
  \multiput(15,85)(60,0){2}{\put(0,0){\cell}}
  \multiput(155,65)(60,0){3}{\put(0,0){\cell}}
  \multiput(5,35)(60,0){3}{\put(0,0){\cell}}
  \multiput(15,35)(60,0){2}{\put(0,0){\cell}}
  \multiput(215,85)(60,0){2}{\put(0,0){\cell}}
  \put(265,85){\cell}
  \put(95,85){\cell}
  \put(95,25){\cell}
  \put(155,25){\cell}
  \put(205,25){\cell}
  \put(25,15){\cell}
  \put(25,35){\cell}
  \put(85,15){\cell}
  \put(145,15){\cell}
  \put(145,25){\cell}
  \put(205,15){\cell}
  \put(195,25){\cell}
  \put(215,25){\cell}
  \put(265,25){\cell}
  \put(255,25){\cell}
  \put(275,25){\cell}
  \put(275,5){\cell}
  \multiput(0,0)(60,0){5}{
  \put(5,45){\line(1,0){40}}
  \put(5,35){\line(1,0){40}}
  \put(15,25){\line(1,0){30}}
  \put(25,15){\line(1,0){20}}
  \put(35,5){\line(1,0){10}}
  \put(5,35){\line(0,1){10}}
  \put(15,25){\line(0,1){20}}
  \put(25,15){\line(0,1){30}}
  \put(35,5){\line(0,1){40}}
  \put(45,5){\line(0,1){40}}
  }
  \multiput(0,60)(60,0){5}{
  \put(5,45){\line(1,0){40}}
  \put(5,35){\line(1,0){40}}
  \put(15,25){\line(1,0){30}}
  \put(25,15){\line(1,0){20}}
  \put(35,5){\line(1,0){10}}
  \put(5,35){\line(0,1){10}}
  \put(15,25){\line(0,1){20}}
  \put(25,15){\line(0,1){30}}
  \put(35,5){\line(0,1){40}}
  \put(45,5){\line(0,1){40}}
  }
  %
  \end{picture}
\end{center}
The five members in the first row
form the subset $\mathcal{E}_{\mu}^{I\!I}(\lambda).$ 
The corresponding elements in $W^{P}$ are
as follows: 
for type $D_{5}:$ $v=5\bar{4}\bar{3}\bar{2}\bar{1},\;
w=135\bar{4}\bar{2},$
for type $C_{4}$ or $B_{4}$:
$v=\bar{4}\bar{3}\bar{2}\bar{1},\;
w=13\bar{4}\bar{2}.$
\end{Exam}

\begin{thm}\label{BCD}
Let $w\leq v\in W^{P},$ and $\lambda\leq \mu$  the
corresponding strict partitions.  
We have the following formulas:
\begin{enumerate} 
\item {\bf Type} $C_n$:
$
[X_w]|_{v}
=\sum_{C\in \mathcal{E}_\mu^I(\lambda)}
\prod_{(i,j)\in C}(\ep_{v(n+j)}-\ep_{v(n-i+1)}),$
\item {\bf Type} $B_n$:
$
[X_w]|_{v}=\sum_{C\in \mathcal{E}_\mu^I(\lambda)}
\prod_{(i,j)\in C}2^{-\delta_{ij}}(\ep_{v(n+j)}-\ep_{v(n-i+1)}),
$
\item {\bf Type} $D_{n+1}$:
$
[X_w]|_{v}
=\sum_{C\in \mathcal{E}_\mu^{I\!I}(\lambda)}
\prod_{(i,j)\in C}(\ep_{v(n+j+2)}-\ep_{v(n-i+2)}).
$
\end{enumerate}
\end{thm}

\begin{Exam} {\rm Let $n=4$ and $\mu=\rho_4.$
We arrange positive roots in $D'_{\mu}$ as follows.}
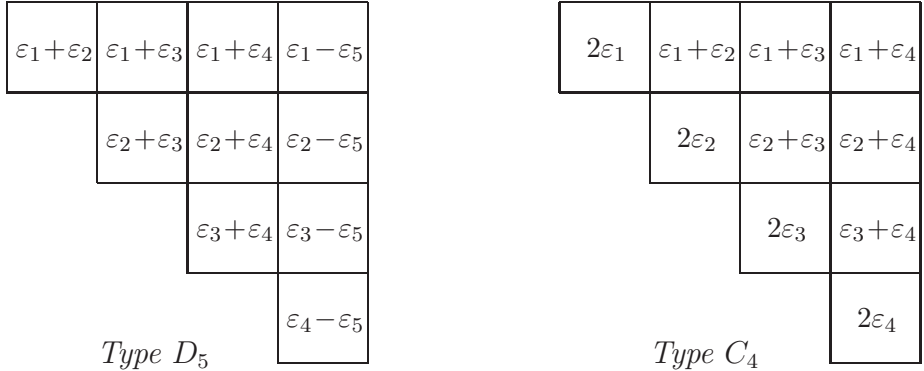
\setlength{\unitlength}{1.2mm}
\begin{center}
  \begin{picture}(60,60)
  \put(5,55){\line(1,0){40}}
  \put(5,45){\line(1,0){40}}
  \put(15,35){\line(1,0){30}}
  \put(25,25){\line(1,0){20}}
  \put(35,15){\line(1,0){10}}
  \put(5,55){\line(0,-1){10}}
  \put(15,55){\line(0,-1){20}}
  \put(25,55){\line(0,-1){30}}
  \put(35,55){\line(0,-1){40}}
  \put(45,55){\line(0,-1){40}}
  \put(6,49){\small{$\ep_{1}\!+\!\ep_2$}}
  \put(16,49){\small{$\ep_{1}\!+\!\ep_3$}}
  \put(26,49){\small{$\ep_{1}\!+\!\ep_4$}}
  \put(36,49){\small{$\ep_{1}\!-\!\ep_5$}}
  \put(16,39){\small{$\ep_{2}\!+\!\ep_3$}}
  \put(26,39){\small{$\ep_{2}\!+\!\ep_4$}}
  \put(36,39){\small{$\ep_{2}\!-\!\ep_5$}}
  \put(26,29){\small{$\ep_{3}\!+\!\ep_4$}}
  \put(36,29){\small{$\ep_{3}\!-\!\ep_5$}}
  \put(36,19){\small{$\ep_{4}\!-\!\ep_5$}}
  \put(15,15){Type $D_5$}
  \end{picture}
  \begin{picture}(60,60)
  \put(5,55){\line(1,0){40}}
  \put(5,45){\line(1,0){40}}
  \put(15,35){\line(1,0){30}}
  \put(25,25){\line(1,0){20}}
  \put(35,15){\line(1,0){10}}
  \put(5,55){\line(0,-1){10}}
  \put(15,55){\line(0,-1){20}}
  \put(25,55){\line(0,-1){30}}
  \put(35,55){\line(0,-1){40}}
  \put(45,55){\line(0,-1){40}}
  \put(8,49){\small{$2\ep_{1}$}}
  \put(16,49){\small{$\ep_{1}\!+\!\ep_2$}}
  \put(26,49){\small{$\ep_{1}\!+\!\ep_3$}}
  \put(36,49){\small{$\ep_{1}\!+\!\ep_4$}}
  \put(18,39){\small{$2\ep_{2}$}}
  \put(26,39){\small{$\ep_{2}\!+\!\ep_3$}}
  \put(36,39){\small{$\ep_{2}\!+\!\ep_4$}}
  \put(28,29){\small{$2\ep_{3}$}}
  \put(36,29){\small{$\ep_{3}\!+\!\ep_4$}}
  \put(38,19){\small{$2\ep_{4}$}}
  \put(15,15){Type $C_4$}
  \end{picture}
\end{center}

{\rm Let $\lambda=(2).$ Then the set $\mathcal{E}_\mu^{I}(\lambda)$
consists of the following six elements:}
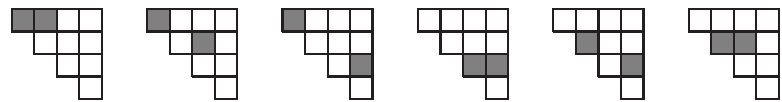
\setlength{\unitlength}{0.3mm}
\begin{center}
  \begin{picture}(300,45) 
  \put(5,35){\cell}
  \put(15,35){\cell}
  \put(65,35){\cell}
  \put(85,25){\cell}
  \put(125,35){\cell}
  \put(155,15){\cell}
  \put(205,15){\cell}
  \put(215,15){\cell}
  \put(255,25){\cell}
  \put(275,15){\cell}
  \put(315,25){\cell}
  \put(325,25){\cell}
  \multiput(0,0)(60,0){6}{
  \put(5,45){\line(1,0){40}}
  \put(5,35){\line(1,0){40}}
  \put(15,25){\line(1,0){30}}
  \put(25,15){\line(1,0){20}}
  \put(35,5){\line(1,0){10}}
  \put(5,35){\line(0,1){10}}
  \put(15,25){\line(0,1){20}}
  \put(25,15){\line(0,1){30}}
  \put(35,5){\line(0,1){40}}
  \put(45,5){\line(0,1){40}}
  }
  \end{picture}
\end{center}
{\rm The first four elements
form the subset $\mathcal{E}_\mu^{I\!I}(\lambda).$}
{\rm In type $D_5$ case we have}
$$[X_w]|_v=
(\ep_1+\ep_2)(\ep_1+\ep_3)
+(\ep_1+\ep_2)(\ep_2+\ep_4)
+(\ep_1+\ep_2)(\ep_3-\ep_5)
+(\ep_3+\ep_4)(\ep_3-\ep_5),
$$
{\rm in type $C_4$ case we have}
$$
[X_w]|_v=2\ep_1(\ep_1+\ep_2)
+2\ep_1(\ep_2+\ep_3)
+2\ep_1(\ep_3+\ep_4)
+2\ep_3(\ep_3+\ep_4)
+2\ep_2(\ep_3+\ep_4)
+2\ep_2(\ep_2+\ep_3).
$$
\end{Exam}

\section{Proof of Theorem \ref{BCD}}\label{PfBCD}
\setcounter{equation}{0}
This section is devoted to the proof of Theorem \ref{BCD}.
Our strategy is the same as in Section \ref{PfEYD}.

We shall give a notion of {\it row-reading expression} for each element
of $W^{P}.$
First we fill in the diagram $D'_{\rho_{n}}$ the simple 
reflections.

{\bf Types} $B_n,C_n$: 
We define a map
$s: D'_{\rho_{n}}\rightarrow \{s_1,\ldots,s_{n}\}$
by 
$$
s(i,j)=\begin{cases}
s_{n} \quad (i=j)\\
s_{n+i-j} \quad (i< j)
\end{cases},
$$
where $1\leq i<j\leq n.$

{\bf Type} $D_{n+1}$: We define a map
$s: D'_{\rho_{n}}\rightarrow \{s_1,\ldots,s_{n+1}\}$
by 
$$
s(i,j)=\begin{cases}
s_{n+1} \quad (i=j, \; i: odd)\\
s_{n} \quad (i=j,\; i: even)\\
s_{n+i-j} \quad (i< j)
\end{cases},
$$
where $1\leq i<j\leq n.$

\begin{Exam} {\rm Let $n=4.$ The map $s$ is illustrated as follows:}
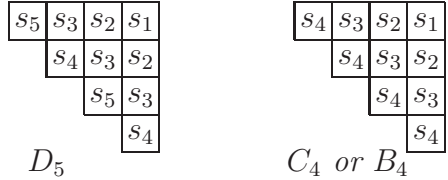
\setlength{\unitlength}{0.5mm}
\begin{center}
  \begin{picture}(50,50)
  \put(5,45){\line(1,0){40}}
  \put(5,35){\line(1,0){40}}
  \put(15,25){\line(1,0){30}}
  \put(25,15){\line(1,0){20}}
  \put(35,5){\line(1,0){10}}
  \put(5,35){\line(0,1){10}}
  \put(15,25){\line(0,1){20}}
  \put(25,15){\line(0,1){30}}
  \put(35,5){\line(0,1){40}}
  \put(45,5){\line(0,1){40}}
  \put(36.5,8.5){$s_4$}
  \put(26.5,18.5){$s_5$}
  \put(16.5,28.5){$s_4$}
  \put(6.5,38.5){$s_5$}
  \put(16.5,38.5){$s_3$}
  \put(26.5,28.5){$s_3$}
  \put(36.5,18.5){$s_3$}
  \put(26.5,38.5){$s_2$}
  \put(36.5,28.5){$s_2$}
  \put(36.5,38.5){$s_1$}
  \put(10,0){$D_5$}
  \end{picture}
  \hspace{1cm}
  \begin{picture}(50,50)
  \put(5,45){\line(1,0){40}}
  \put(5,35){\line(1,0){40}}
  \put(15,25){\line(1,0){30}}
  \put(25,15){\line(1,0){20}}
  \put(35,5){\line(1,0){10}}
  \put(5,35){\line(0,1){10}}
  \put(15,25){\line(0,1){20}}
  \put(25,15){\line(0,1){30}}
  \put(35,5){\line(0,1){40}}
  \put(45,5){\line(0,1){40}}
  \put(36.5,8.5){$s_4$}
  \put(26.5,18.5){$s_4$}
  \put(16.5,28.5){$s_4$}
  \put(6.5,38.5){$s_4$}
  \put(16.5,38.5){$s_3$}
  \put(26.5,28.5){$s_3$}
  \put(36.5,18.5){$s_3$}
  \put(26.5,38.5){$s_2$}
  \put(36.5,28.5){$s_2$}
  \put(36.5,38.5){$s_1$}
  \put(3,0){$C_4$ or $B_4$}
  \end{picture}
\end{center}
\end{Exam}

Let $\lambda$ be the strict partition 
associated to $w\in W^P.$
We define the row-reading map $\varphi:D'_\lambda\rightarrow \{1,2,\ldots,k\}$
as in Section \ref{EYD}.
We define the word
$$
s_{i_1}\cdots s_{i_{k}},
$$
where $s_{i_j}=s(\varphi^{-1}(j)).$ 
Then we have 
$w=s_{i_1}\cdots s_{i_{k}}$,
which is a reduced  expression for $w\in W^{P}.$
For each subset $C$ of $D'_\mu,$
we define an element $w_C\in W$
in the same way in Section \ref{EYD}.
Let $w\leq v\in W^P.$
The set 
$$
\{C\subset D'_\mu\;|\; \sharp C=\ell(w),\;w_C=w\}
$$
is denoted by $\mathcal{R}_v^{I}(w)$ (resp. $\mathcal{R}_v^{I\!I}(w)$ )
for types $C_n,B_n$ ( type $D_{n+1}$). 

The key result to prove Theorem \ref{BCD} is the following.
\begin{prop}\label{key} Let $w\leq v\in W^P$ and $\lambda\leq \mu$ the associated
strict partitions.
Then we have
$\mathcal{E}_\mu^{*}(\lambda)=
\mathcal{R}_v^*(w)
$
for $*=I,I\!I.$
\end{prop}

\begin{lem}\label{shiftedCC'}
Let $C,C'$ be subsets of a shifted Young diagram $D'_\mu$
such that $C'$ is obtained from $C$ by an elementary excitation of type $*=I,I\!I.$
If $C$ belongs to $\mathcal{R}_v^*(w)$
then we also have $C'\in \mathcal{R}_v^*(w).$ 
\end{lem}
{\it Proof.}
Suppose $C'$ is obtained from $C$   
by an elementary excitation of Type $I\!I$
occurring at $(i,i).$ 
Any simple reflection contained 
in the region $R_\flat\cap C$
are one of $s_{1},\ldots,s_{n-2}$
each of which commutes with $s_{n+1}$ and $s_n.$
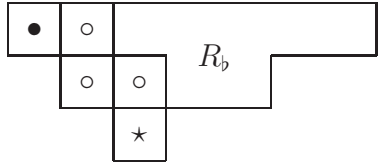
\begin{center}
\setlength{\unitlength}{0.7mm}
\begin{picture}(160,40)
\put(20,30){\line(1,0){70}}
\put(30,10){\line(1,0){40}}
\put(40,0){\line(1,0){10}}
\put(20,20){\line(1,0){30}}
\put(70,10){\line(0,1){10}}
\put(90,20){\line(0,1){10}}
\put(70,20){\line(1,0){20}}
\put(30,10){\line(0,1){20}}
\put(40,0){\line(0,1){30}}
\put(50,0){\line(0,1){20}}
\put(20,20){\line(0,1){10}}
\put(43.5,3.5){$\star$}
\put(23.5,23.5){$\bullet$}
\put(33.5,23.5){$\circ$}
\put(33.5,13.5){$\circ$}
\put(43.5,13.5){$\circ$}
\put(56,18){$R_\flat$}
\end{picture}
\end{center}
The positions indicated by $\circ$ 
are vacant by the definition of 
elementary excitation.
Therefore the word $w_{C'}$ can obtained from $w_{C'}$  
by using only commuting
relations.  
Hence we have $C'\in \mathcal{R}_v(w).$
The arguments for the other cases are simpler 
than this and we leave them to the reader.
\hfill $\square$

\begin{lem}\label{shiftedRE}
Let $C\in \mathcal{R}_v^*(w)$ be such that $D_\lambda'\ne C.$
Then there exists $C'\in\mathcal{R}_v^*(w)$ such that 
$C$ is obtained from $C'$ by excitation of type $*.$
\end{lem}
{\it Proof.} Similar to the proof of Lemma \ref{RE}.
\hfill $\square$

\bigskip

{\it Proof of Proposition \ref{key}.}
This is immediate from Lemmas \ref{shiftedCC'}, \ref{shiftedRE}.
\hfill $\square$

\bigskip

Now we calculate the $\beta$-sequence in 
Proposition \ref{AJS_B}.

\begin{lem} \label{betaBCD}
Let $\mu\in \mathcal{SP}_{\rho_n}$
and $v$ the associated element of $W^P.$
Denote by 
$\varphi:\mu\longrightarrow\{1,2,\ldots,k\} $
the row-reading map of $\mu.$
For $(i,j)\in D'_{\mu}$ set $\beta_{i,j}=\beta_{\varphi(i,j)},$
where $\beta_t\:(1\leq t\leq l)$ are defined
in (\ref{beta}). 
We have the following formulas:
\begin{enumerate}
\item {\bf Type} $C_{n}$:
$\beta_{i,j}=\ep_{v(n+j)}-\ep_{v(n-i+1)},$
\item {\bf Type} $B_{n}$:
$\beta_{i,j}=2^{-\delta_{ij}}
(\ep_{v(n+j)}-\ep_{v(n-i+1)}),$
\item {\bf Type} $D_{n+1}$:
$\beta_{i,j}=\ep_{v(n+2+j)}-\ep_{v(n-i+2)}.$
\end{enumerate}
\end{lem}
{\it Proof.}
Similar to the proof of
Lemma \ref{betalem1} and \ref{corner}. \hfill $\square$

\bigskip

{\it Proof of Theorem \ref{BCD}.}
The same as the proof of Theorem \ref{EYD}.
\hfill $\square$

\section{Factorial $Q$- and $P$-functions}\label{SectQP}
\setcounter{equation}{0}
In this section, we first define the factorial
$Q$- and $P$-functions
and present some fundamental properties.
Then we use them to express 
the restriction of the equivariant Schubert classes
for the Lagrangian and orthogonal Grassmannians.
As an application we give the equivariant Giambelli formula.

\subsection{Factorial $Q$- and $P$-functions}
We first recall the definition of 
factorial $Q$- and $P$-functions due to 
Ivanov \cite{Iv}.
Let $\mathcal{SP}(n)$ denote the 
set of strict partitions
$\lambda=(\lambda_1>\cdots>\lambda_r>0)$
with $r\leq n.$
\begin{Def} Let $\lambda\in \mathcal{SP}(n).$
Put
\beqn
P^{(n)}_\lambda(x|a)=
\frac{1}{(n-r)!}
\sum_{w\in S_n}w\left(
(x_1|a)^{\lambda_1}\cdots
(x_r|a)^{\lambda_r}
\prod_{1\leq i\leq r,\;i\leq j\leq n}
\frac{x_i+x_j}{x_i-x_j}
\right),\label{defP}
\eeqn
where $w\in S_n$ acts as a permutation of 
variables $x_1,\ldots,x_n.$
We also put $Q^{(n)}_\lambda(x|a)=2^{r}P^{(n)}_\lambda(x|a).$
\end{Def}
The rational expression in the 
right hand side of (\ref{defP}) is 
actually a polynomial in $x_1,\ldots,x_n$
and $a_1,a_2,\ldots,a_{\lambda_1},$
homogeneous of degree $|\lambda|.$
In particular we have
\beqn
P^{(n)}_1(x|a)=\begin{cases}
x_1+\cdots+x_n & \mbox{if}\;n\;\mbox{is even}\\
x_1+\cdots+x_n-a_1& \mbox{if}\;n\;\mbox{is odd}
\end{cases}.\label{P1}
\eeqn

In \cite{Iv}, the first parameter $a_1$ is assumed to be zero in the most 
part of the argument. 
However, we need $a_1$ for the later use.
In order to generalize Ivanov's results to the case of 
non-zero $a_1$ the following is fundamental
(cf. \cite{Iv}, Proposition 2.6).

\begin{prop}\label{Stab}(Stability mod $2$) 
For any $\lambda\in \mathcal{SP}(n),$
we have
\beqn
P^{(n+2)}_\lambda(x_1,\ldots,x_n,0,0|a)
=
P^{(n)}_\lambda(x_1,\ldots,x_n|a).\label{stab2}
\eeqn
\end{prop}
{\it Proof.} We set
$$
F_n=(x_1|a)^{\lambda_1}\cdots
(x_r|a)^{\lambda_r}
\prod_{1\leq i\leq r,\;i\leq j\leq n}
\frac{x_i+x_j}{x_i-x_j}.
$$
Let $w\in S_{n+2}$, and consider the term $w(F_{n+2}).$
If $w(i)\ne n+1,n+2$ for all $i=1,\ldots,r$, then
the rational function $w(F_{n+2})$ has 
no pole along the hyperplane $x_{n+1}-x_{n+2}=0.$
So we can substitute $x_{n+1}=x_{n+2}=0$ to 
such a rational function. As the result of 
the substitution we have $w'(F_n)$ for some $w'\in S_n.$
In fact, $w'$ is obtained by eliminating 
$n+1,n+2$ from the sequence
$w(1),\ldots,w(n+2).$
By this correspondence, there are $(n-r+2)(n-r+1)$ elements 
of $S_{n+2}$ giving the same $w'\in S_n.$
Thus our task is to 
show that all the 
remaining terms cancel out.
 
Suppose $i=w^{-1}(n+1)\leq r$ or $j=w^{-1}(n+1)\leq r$ and
set $w'=w\cdot (i,j).$
Then both $w(F_{n+2})$ and $w'(F_{n+2})$
has simple pole along $x_{n+1}-x_{n+2}=0,$ however, the sum
$w(F_{n+2})+w'(F_{n+2})$ changes sign 
when we permute the variables $x_{n+1}$ and $x_{n+2}$, and 
hence regular along $x_{n+1}-x_{n+2}=0.$
Thus we can substitute $x_{n+1}=x_{n+2}=0$ to the sum and 
obtain zero because it has a factor $x_{n+1}+x_{n+2}.$
\hfill $\square$

\subsection{Vanishing property}
Here we present a vanishing property
of $P^{(n)}_\lambda(x|a).$ 
For strict partition $\lambda=(\lambda_1>\cdots>\lambda_r>0)$ 
with $r\leq n$ 
define an $n$-tuple
$$
a_\lambda=\begin{cases}
(a_{\lambda_1+1},\ldots,a_{\lambda_r+1},0,\ldots,0)& \mbox{if}\quad n-r\quad
\mbox{is even}\\
(a_{\lambda_1+1},\ldots,a_{\lambda_r+1},a_1,0,\ldots,0)& \mbox{if}\quad n-r\quad
\mbox{is odd}
\end{cases}.
$$

\begin{prop}\label{vanishP} We have
$P^{(n)}_\lambda(a_\mu|a)=0$ unless $\mu\geq \lambda.$
\end{prop}
{\it Proof.} By Proposition \ref{Stab}, we may
assume $a_\mu$ has no zero entries.
Then the proposition follows from definition.
\hfill $\square$

\bigskip

We record a formula for 
$P^{(n)}_\lambda(a_\lambda|a)$ below,
although
this result is not used  
in this paper.
Let $\{\mu_1,\ldots,\mu_{n-r}\}$ be
the subset of $\{2,3,\ldots,n+1\}$
complementary to $\{\lambda_i+1\;|\;1\leq i\leq r\},$
where we arrange the elements 
in increasing order, i.e. $\mu_1<\cdots<\mu_{n-r}.$
Define $\mu_0=1$ if $n-r$ is odd
and $\mu_0=\bar{1}$ if $n-r$ is even.
Then we denote the sequence
$$
\lambda_1+1,\ldots,\lambda_r+1,\mu_0,\overline{{\mu}_1},\ldots,
\overline{{\mu}_{n-r}},
$$
by $\kappa_1,\ldots,\kappa_{n+1}.$
Put
$$
H_\lambda(a)=
\prod_{(i,j)\in D'_\lambda}(a_{\kappa_i}+a_{\kappa_{j+1}}),
$$
where we set $a_{\bar{i}}=-a_i.$

\begin{prop} For a strict partition $\lambda\in \mathcal{SP}(n)$ we have
$P^{(n)}_\lambda(a_\lambda|a)=H_\lambda(a).$
\end{prop}

\subsection{Pieri's rule and Pfaffian formulas}
We collect here some basic facts on the factorial
$P$-Schur functions,
which we shall make use of
in the next subsection.
The proof for $a_1=0$ case is given in \cite{Iv}.
The same proof works for the general $a_1$ case
using the vanishing property. 

\begin{prop}\label{PieriP} (\cite{Iv}) 
For a strict partition $\lambda\in \mathcal{SP}(n)$ 
we have
$$
\left(P_{1}^{(n)}(x|a)-P^{(n)}_1(a_\lambda|a)\right)
P_\lambda^{(n)}(x|a)=\sum_{\lambda\to\nu}P_\nu^{(n)}(x|a),
$$
where the sum is over $\nu\in \mathcal{SP}(n)$ such 
that $\nu\geq \lambda$ and $|\nu|=|\lambda|+1.$
\end{prop}

For any
$\lambda=(\lambda_1>\cdots>\lambda_r>0)\in \mathcal{SP}(n),$
set $r_0(\lambda)=r$ if $r$ is even and $r_0(\lambda)=r+1$
if $r$ is odd, and then we put $\lambda_{r+1}=0.$
\begin{lem}\label{Pf}(\cite{Iv}) 
For a strict partition $\lambda\in \mathcal{SP}(n)$ 
we have
$$
P^{(n)}_\lambda(x|a)=\mathrm{Pf}
(P^{(n)}_{\lambda_i,\lambda_j}(x|a))_{1\leq i<j\leq r_0(\lambda)}.
$$
\end{lem}

\subsection{Closed formula for $[X_w]|_{v}$ and equivariant Giambelli formula}
For each $v\in W^P$, we shall define $n$-tuple $x_v$
in the following way: 

{\bf Types} $C_n$ and $B_n$:
Let $v\in W^{P},$ and 
$$
\{i\;|\; 1\leq i\leq n,\;
v(i)>n\}
=\{\bar{j}_1,\ldots,\bar{j}_k\},
$$
where $1\leq j_1<\cdots<j_k\leq n.$
Then we put
$$
x_v=(\ep_{i_1},\ldots,\ep_{j_k},0,\ldots,0).
$$

{\bf Type} $D_{n+1}$: Let $v\in W^{P},$ and
$$
\{i\;|\; 1\leq i\leq n+1,\;v(i)>n+1\}
=\{\bar{j}_1,\ldots,\bar{j}_k\},
$$
where $1\leq j_1<\cdots<j_k\leq n+1.$
Note that $k$ is even.
If $n$ is even we put
$$
x_v=(\ep_{j_1},\ldots,\ep_{j_k},0,\ldots,0).
$$
If $n$ is odd, let $r=r(\lambda)$, where $\lambda$ 
is the strict partition corresponding to $v$,
and put 
$$
x_v=\begin{cases}
(\ep_{j_1},\ldots,\ep_{j_r},0,\ldots,0) \quad 
&\mbox{if}\; n-r \;\mbox{is odd}\\
(\ep_{j_1},\ldots,\ep_{j_r},-\ep_{n+1},0,\ldots,0) \quad 
&\mbox{if}\; n-r \;\mbox{is even}
\end{cases}.
$$

\begin{lem}
Let $v\in W^P.$ We have
\begin{enumerate}
\item
{\bf Type} $C_n$:
$
[X_{s_n}]|_{v}
=Q^{(n)}_1(x_v|0,\ep_n,\ldots,\ep_2),
$
\item {\bf Type} $B_n$:
$
[X_{s_n}]|_{v}
=P^{(n)}_1(x_v|0,\ep_n,\ldots,\ep_2),
$
\item {\bf Type} $D_{n+1}$:
$
[X_{s_{n+1}}]|_v=P^{(n)}_1(x_v|(-1)^n\ep_{n+1},\ep_n,\ldots,\ep_2).
$
\end{enumerate}
\end{lem}

{\it Proof.} Consider Type $D_{n+1}.$
By Proposition \ref{xi}, (4), and (\ref{pi}),
we have
$$
[X_{s_{n+1}}]|_v=\varpi_{n+1}-v(\varpi_{n+1})
=\sum_{1\leq i\leq n+1,\,v(i)>n+1}\ep_{i}.
$$
Now recall the form of $P_1(x|a)$ is given by (\ref{P1}).
Then it is easy to check our 
formula.
Types $C_n$ and $B_n$ are similar and much simpler.
\hfill $\square$

\bigskip

Now we can state the main result of this 
section.
\begin{thm}\label{BCD-closed}
Let $w\leq v\in W^{P},$ and $\lambda\leq \mu$  the
corresponding strict partitions. 
We have the following formulas:
\begin{enumerate}
\item
{\bf Type} $C_n$ (\cite{Ik}):
$
[X_w]|_{v}
=Q^{(n)}_\lambda(x_v|0,\ep_n,\ldots,\ep_2),
$
\item {\bf Type} $B_n$:
$
[X_w]|_{v}
=P^{(n)}_\lambda(x_v|0,\ep_n,\ldots,\ep_2),
$
\item {\bf Type} $D_{n+1}$:
$
[X_w]|_{v}
=P^{(n)}_\lambda(x_v|(-1)^n\ep_{n+1},\ep_n,\ldots,\ep_2).
$
\end{enumerate}
\end{thm}

{\it Proof.} The result for Type $C_n$ has been proved in \cite{Ik}.
We consider Type $D_{n+1}.$
Assume $n$ is even for simplicity.
The odd case is left for the reader. 
Fix $\mu\in \mathcal{SP}_{\!\rho_n}.$ 
Consider the following system of equations for
the functions 
$F_\lambda(\ep_1,\ldots,\ep_{n+1})\;(\lambda\in \mathcal{SP}_{\mu}):$
\beqn
\left(\sum_{i=1}^{r_0(\mu)}\ep_{n-\mu_i+1}
-\sum_{i=1}^{r_0(\lambda)}\ep_{n-\lambda_i+1}\right)F_\lambda(\ep)
=\sum_{\nu}F_\nu(\ep)\quad
\mbox{for all}\; \lambda\in \mathcal{SP}_\mu\setminus \{\mu\},
\label{recBCD}
\eeqn
where the summation in the right hand side runs through
$\nu\in \mathcal{SP}_{\mu}$ such that $\nu\geq \lambda$ and 
$|\nu|=|\lambda|+1.$ 
By the equations (\ref{rec}) together with 
the initial condition $F_\phi(\ep_1,\ldots,\ep_n)=1,$
the set of functions $F_\lambda(\ep)\;(\lambda\in \mathcal{SP}_\mu)$ 
is characterized.

By equation (\ref{ParaXi}) for type $D_{n+1}$,
$[X_w]|_v=\xi^{w}(v)$ (cf. Proposition \ref{A}) 
satisfy equation (\ref{recBCD}) as well as the initial condition
(cf. Proposition \ref{xi}, (1)).
Note that the coefficients $\langle w(\varpi_{n+1}),\beta^\vee\rangle$,
the Chevalley multiplicity, 
is equal to one 
in the right hand side of equation (\ref{ParaXi}),
i.e., $G/P$ in this case is `minuscule'.
On the other hand, the functions on the right hand side of
the formula in the theorem
also satisfy (\ref{recBCD}).
This fact is a consequence of Prpoposition \ref{PieriP} and 
the vanishing property (Proposition \ref{vanishP}) of $P^{(n)}_\lambda(x|a).$
The initial condition is also satisfied.
Therefore we have the theorem.
The type $B_n$ case is quite similar to the case of type $D_{n+1}.$
\hfill $\square$

\bigskip

The following result is a direct 
consequence of Theorem \ref{BCD-closed}
and Lemma \ref{Pf}.

\begin{cor}\label{Giam}(Equivariant Giambelli) 
Let $G/P$ of type $C_n,B_n,D_{n+1}.$ 
For any $\lambda\in \mathcal{SP}_{\!\rho_n},$
we have
$$
[X_\lambda]
=\mathrm{Pf}\left([X_{\lambda_i,\lambda_j}]\right)_{1\leq i<j\leq r_0(\lambda)},
$$
where we denote by $X_\lambda$ the Schubert 
variety corresponding to $\lambda.$
\end{cor}

{\it Proof.}
The proof is the same as that given in \cite{Ik}.
\hfill $\square$

\section{Multiplicity of a singular point in a Schubert variety}\label{SectMult}
\setcounter{equation}{0}
Another application is to the multiplicity of a singular point in
a Schubert variety.
We denote the multiplicity of the variety 
$X_w$ at $e_v$ by $m_v(X_w).$
We will explain the relationship between 
$[X_w]|_v$ and $m_v(X_v).$
Then we discuss some implication
of our result on $m_v(X_v).$

Let $R_{uni}(P)$ be the unipotent radical
of $P$,
and let $R_P^{+}$ the subset of $R^{+}$ 
defined by $R^{+}_P=\{\beta\in R\;|\;U_\beta\subset R_{uni}(P)\},$
where $U_\beta$ is the root subgroup associated to $\beta.$
Given $v\in W^P.$ Let $U^{-}_v$ be the subgroup of $G$ generated by the
root subgroups $U_{-\beta},\;\beta\in v(R^{+}\setminus R_P^{+}).$
Under the map $G\rightarrow G/P,$ $U^{-}_v$ is 
mapped isomorphically onto its image $\mathcal{U}_v:=U_{v}^{-}e_v$ which
is a canonical $T$-stable affine neighborhood of $e_v$ 
with a coordinate system $\{x_{-\beta}\;|\;
\beta\in v(R^{+}\setminus R_P^{+})\}.$

Let $G/P$ be either of types $A_{n-1},\,C_{n},\,D_{n}.$
For $w,v\in W^P,$ $v\geq w$, let us denote 
$Y_{w,v}=X_w\cap \,\mathcal{U}_{v}.$
It is known that in these cases the defining 
ideal of the affine variety
$Y_{w,v}$ is homogeneous in our coordinate system,
i.e. $Y_{w,v}$ is a cone in $\mathcal{U}_{v}.$
%
Actually we have a  
one parameter subgroup 
$\phi_{v}\in \mathrm{Hom}(\C^{\times},T)$ 
that $\phi_{v}(t)\;(t\in \C^{\times})$ 
acts on $\mathcal{U}_{v}$ by dilation.
Then there exists $h_v \in \mathrm{Lie}(T)$
such that
$h_v(\beta)=-1$ for all $\beta\in v(R^{+}\setminus R_P^{+}).$ 
For example, if $v\in W(C_{n})^{P_{n}}$ 
then $h_v$ is given by 
$h_v(\ep_i)=\frac{1}{2}$ if $v(i)>n$ and 
$-\frac{1}{2}$ if $v(i)\leq n$,
for $1\leq i\leq n.$
Note that each $[X_w]|_v\in S$, being a polynomial function
on $\mathrm{Lie}(T)$, can be evaluated at $h_v.$

\begin{prop}\label{special} Let $G/P$ be either of 
types $A_{n-1},\,C_{n},\,D_{n}$ and $w\leq v\in W^P.$
The value of $[X_w]|_v$ evaluated at $h_v$ gives 
the multiplicity of the variety 
$X_w$ at $e_v.$
\end{prop}
{\it Proof.} Consider the neighborhood $\mathcal{U}_v$ of $e_v.$
The multiplicity $m_v(X_w)$ is determined by 
the Poincar\'e series of the coordinate ring $\C[Y_{w,v}],$
which is given 
by restricting the formal character $\mathrm{ch}\,\C[Y_{w,v}]$ as a $T$-module
to the one parameter subgroup 
$\phi_{v}$ associated with $h_v.$ 
Then $[X_w]|_v$ evaluated at $h_v$ is nothing 
but the classical multiplicity 
in the sense of Hilbert-Samuel. 
\hfill$\square$

{\it Remark.} There are
recurrence relations for $m_{v}(X_{w})$
given by Lakshmibai and Weyman \cite{LW},
which can be obtained by specializing 
equations (\ref{rec}), (\ref{recBCD}) 
at $h_{v}.$
So we have another proof of the above 
Proposition.

\bigskip

\begin{cor}\label{multEYD} (\cite{KL},\cite{KR},\cite{Kr1}) 
Let $G/P=\mathcal{G}_{d,n},$ $w\leq v\in W^{P_d},$
and $\lambda\leq \mu\in \mathcal{P}_{d,n}$ be the 
corresponding partitions.
Then we have
\beqn
m_{v}(X_w)=\sharp \mathcal{E}_\mu(\lambda).\label{classical_A}
\eeqn
\end{cor}

{\it Proof.} This is immediate from Proposition \ref{special} 
and Theorem \ref{EYD},
since each summand in the right hand side 
of the formula becomes one, when
specialized at $h_v.$
\hfill$\square$

\bigskip

{\it Remark 2.} This formula can be obtained as a 
direct consequence of the result in \cite{KL}, \cite{KR}
describing the Gr\"obner basis of the defining ideal
of $Y_{w,v},$ 
together with a combinatorial argument in 
\cite{Kr1}.

\begin{cor}\label{MultCD} Let $G/P$ be the Grassmannian of types $C_n$ or $D_{n+1}.$
Let $w\leq v\in W^{P},$
and $\lambda\leq \mu\in \mathcal{SP}_{\!\rho_n}$ be the 
corresponding strict partitions.
Then we have
\begin{enumerate}
\item (\cite{GR},\cite{Kr2}) {\bf Type} $C_{n}$
\beqn
m_{v}(X_w)=\sharp\mathcal{E}_\mu^{I}(\lambda),\label{classical_C}
\eeqn
\item {\bf Type} $D_{n+1}$
\beqn
m_{v}(X_w)=\sharp\mathcal{E}_\mu^{I\!I}(\lambda).\label{classical_D}
\eeqn
\end{enumerate}
\end{cor}

{\it Proof.} The same as Corollary \ref{multEYD}.
\hfill$\square$

\begin{Exam} {\rm (cf. Example \ref{strictEYD}) Let 
$v=5\bar{4}\bar{3}\bar{2}\bar{1},\;
w=135\bar{4}\bar{2}$ (type $D_{5}$), and
$v=\bar{4}\bar{3}\bar{2}\bar{1},\;
w=13\bar{4}\bar{2}$
(type $B_{4}$ or $C_{4}$).
Then we have $m_{v}(X_{w})=10$
for $C_{4}$ and
$m_w(X_v)=5$
for $B_{4},D_{5}.$}
\end{Exam}

{\it Remark 3.} 
We should make a remark on (\ref{classical_C}) similar to 
{\it Remark 2}.
Ghorpade and Raghavan \cite{GR} has given a detailed 
description of the Gr\"obner basis of the defining ideal
of $Y_{w,v}$ for the Lagrangian Grassmannian.
Then (\ref{classical_C}) can be 
derived from a result in \cite{Kr2}.
As for type $D_{n+1},$ formula (\ref{classical_D})
seems to be new.

\bigskip

We close this section by pointing out the following fact.

\begin{prop}
With notations as in
Corollary \ref{MultCD}, we have
$$
m_{v}(X_w)=\mathrm{Pf}(m_{v}(X_{\lambda_i,\lambda_j}))_{1\leq i<j\leq r_0(\lambda)}.
$$
\end{prop} 

{\it Proof.} This is
immediate from Corollary \ref{Giam} and 
Proposition \ref{special}.
\hfill $\square$

\section{Lattice paths method}\label{SectLPM}
\setcounter{equation}{0}
As a supplementary discussion, we will show
a direct combinatorial route from the combinatorial formula (Theorem \ref{BCD})
to the equivariant Giambelli formula (Corollary \ref{Giam}).
Our argument relies on the ``lattice
path method'' due to Stembridge. In order to
deal with the case of type $D_{n+1}$ also, 
we slightly modify the results in \cite{St2}. 

\subsection{Perfect matchings}
Let $r$ be an even positive integer.
Denote by $\mathcal{M}_r$ the set of 
all perfect matchings of the set $\{u_1,\ldots,u_r\}.$
We denote a perfect matching $\sigma\in \mathcal{M}_r$ in
such a way that $\sigma=\{(u_{i_k},u_{j_k})\}_{k=1,\ldots,r/2},$
where $i_k<j_k\;(1\leq k\leq r/2),\;i_1<\cdots<i_{r/2},$
and put $I_\sigma=\{i_1,\ldots,i_r\}.$ 
Let $\sigma_0$ denote the `identical' perfect matching 
$\{(u_{2i-1},u_{2i})\}_{i=1}^r.$
\begin{lem}\label{invo}
There is an involution on $\mathcal{M}_r\setminus \{\sigma_0\}$
$(\sigma\mapsto \sigma^{*})$ such that $\sgn(\sigma)=-\sgn(\sigma^{*})$
and $I_\sigma=I_{\sigma^{*}}.$
\end{lem}
{\it Proof.} Take the smallest number $k$ such that 
$(u_i,u_{2k-1}),(u_j,u_{2k})\in \sigma.$
Define a perfect matching $\sigma'$ by replacing 
$(u_i,u_{2k-1}),(u_j,u_{2k})$ in $\sigma$ 
into
$(u_j,u_{2k-1}),(u_i,u_{2k}).$
Then it is easy to see that $\sigma^{*}$ has the desired property.
\hfill $\square$

\subsection{Modified version of Stembridge's argument}
Let $\mu=(\mu_1,\ldots,\mu_l)\in \mathcal{SP}_{\!\rho_n}.$
We define a directed graph as follows.
The vertex set is
$$
D'_\mu\cup (D'_\mu+\pmb{a}),
$$
where $\pmb{a}=(\frac{1}{2},-\frac{1}{2}).$
Direct an edge 
from $u$ to $v$ if (1) $v-u=\pmb{b},$ 
with $\pmb{b}=(-\frac{1}{2},-\frac{1}{2})$
or (2) $u\in D'_\mu$ and $v-u=\pmb{a}.$
If we take $\mu=\rho_n$, the directed graph is the following.
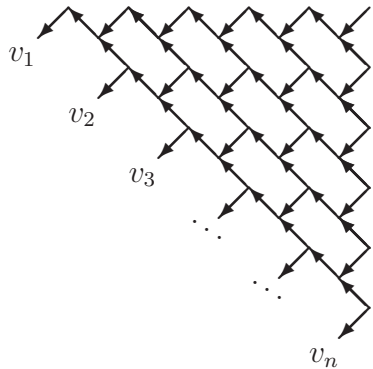
\setlength{\unitlength}{0.8mm}
\begin{center}
  \begin{picture}(90,70)
  \thicklines
  \multiput(75,50)(-10,0){5}{\vector(-1,1){5}}
  \multiput(75,40)(-10,0){4}{\vector(-1,1){5}}
  \multiput(75,30)(-10,0){3}{\vector(-1,1){5}}
  \multiput(75,20)(-10,0){2}{\vector(-1,1){5}}
  \multiput(75,10)(0,10){4}{\vector(-1,1){5}}
  \multiput(70,15)(0,10){4}{\vector(-1,1){5}}
  \multiput(60,25)(0,10){3}{\vector(-1,1){5}}
  \multiput(50,35)(0,10){2}{\vector(-1,1){5}}
  \multiput(40,45)(0,10){2}{\vector(-1,1){5}}
  \multiput(70,55)(-10,0){5}{\vector(-1,1){5}}
  \multiput(75,60)(-10,0){6}{\vector(-1,-1){5}}
  \multiput(75,50)(-10,0){5}{\vector(-1,-1){5}}
  \multiput(75,40)(-10,0){4}{\vector(-1,-1){5}}
  \multiput(75,30)(-10,0){3}{\vector(-1,-1){5}}
  \multiput(75,20)(-10,0){2}{\vector(-1,-1){5}}
  \multiput(75,10)(-10,0){1}{\vector(-1,-1){5}}
  \put(15,51){$v_1$}
  \put(25,41){$v_2$}
  \put(35,31){$v_3$}
  \put(45,21){$\ddots$}
  \put(55,11){$\ddots$}
  \put(65,1){$v_n$}
  \end{picture}
\end{center}
Given a strict partition $\lambda=(\lambda_1,\ldots,\lambda_r)$
such that 
$\lambda\leq \mu$,
let $(u_1,\ldots,u_r)$ be the $r$-tuple of 
vertices defined  
as follows.
If $\mu=\rho_n$ then we put
$$
u_i=(n+1-\lambda_i,n).
$$ 
In general, we put $u_i=(i+m_i,\lambda_i+i-1+m_i)$
where $m_i$ is the largest 
non-negative integer such that $(i+m_i,\lambda_i+i-1+m_i)\in D'_{\mu}.$
Define
$v_i=(i,i)+\pmb{a}$
for $1\leq i\leq  l$ and set $J=\{v_1,\ldots,v_r\},$ and also
$$
J_1=\{v_i\;|\; i\;\mbox{is odd}\},\quad
J_2=\{v_i\;|\; i\;\mbox{is even}\}.
$$
For each subset $I$ of $J$ we denote 
by $\mathcal{P}(u_i,I)$ the set of paths
starting from $u_i$ to any vertex in $J.$
Let $\mathcal{P}_\mu^{I}(\lambda)$ denote the set of 
$r$-tuple of paths $(p_1,\ldots,p_r)$ such that 
$p_i\in \mathcal{P}(u_i,J)$ for $1\leq i \leq r.$
Let $\mathcal{P}_\mu^{I\!I}(\lambda)$ denote the subset of 
$\mathcal{P}_\mu^{I}(\lambda)$
such that $$
p_i\in \mathcal{P}(u_i,J_{i})\quad \mbox{for all}\;1\leq i\leq r,
$$
where the index $i$ in $J_i$ is viewed as modulo $2.$
For $*=I,I\!I,$
we define $\mathcal{N}_\mu^*(\lambda)$
to be the set of 
$r$-tuple $\pmb{p}=(p_1,\ldots,p_r)$ of 
non-intersecting paths
in $\mathcal{P}_\mu^*(\lambda).$

For $u=(i,j)\in D'_\mu$, we assign an arbitrary weight $a_{ij}$ 
(any element in a fixed commutative ring) to 
the edge $u\to u+\pmb{a}.$
All the other edge
is assigned the weight $1.$
For $r$-tuple of paths $\pmb{p}=(p_1,\ldots,p_r)$ 
in $\mathcal{P}_\mu^*(\lambda),$
let $w(p_i)$ denote the product of all weights
which $p_i$ go through and put $w(\pmb{p})=w(p_1)\cdots w(p_r).$
Let 
$$
\mathcal{W}_\mu^*(\lambda)=\sum_{\pmb{p}\in\mathcal{N}_\mu^*(\lambda)}
w(\pmb{p})
$$
denote the corresponding generating function.

\begin{prop}\label{LPMPf}(\cite{St2}) 
Let $\lambda\leq \mu\in \mathcal{SP}_{\!\rho_{n}}.$
We have
\beqn
\mathcal{W}_\mu^*(\lambda)=\mathrm{Pf}\;\mathcal{W}_\mu^*
(\lambda_i,\lambda_j)_{1\leq i<j\leq r_{0}(\lambda)}.\label{WPfaff}
\eeqn
\end{prop}

{\it Proof.} For the case of $*=I$, 
we can apply Theorem 3.1 in \cite{St2} directly.
Here we consider the case of $*=I\!I.$
We may assume that $r=r_{0}(\lambda)$ is even
(see {\it Remark\/} after Theorem 3.1, \cite{St2} ).
Given a perfect matching $\sigma$ of $\{u_1,\ldots,u_r\}.$
An $r$-tuple of paths $\pmb{p}=(p_1,\ldots,p_r)
\in \prod_{i=1}^r\mathcal{P}(u_i,J)$ is 
$\sigma$-{\it admissible\/} if
for each $(u_i,u_j)\in \sigma$ with $i<j$,
$p_i$ and $p_j$ do not intersect and
$p_i\in \mathcal{P}(u_i,J_1),$ 
$p_j\in \mathcal{P}(u_j,J_2).$
Denote the set
$
\hat{\mathcal{P}}
$
of pairs $(\sigma,\pmb{p}),$ where
$\sigma$ is a perfect matching of $\{u_1,\ldots,u_r\}$,
and $\pmb{p}=(p_1,\ldots,p_r)$ is an $r$-tuple of paths 
that is $\sigma$-admissible.
If we assign the weight $\mathrm{sgn}(\sigma)w(\pmb{p})$
to $(\sigma,\pmb{p})\in \hat{\mathcal{P}},$
then the Pfaffian of the right hand side
of (\ref{WPfaff})
is equal to the following generating function:
$$
G(\hat{\mathcal{P}})=\sum_{(\sigma,\,\pmb{p})\in \hat{\mathcal{P}}}
\mathrm{sgn}(\sigma)w(\pmb{p})
$$
for $\hat{\mathcal{P}}.$
Let $\hat{\mathcal{P}}^{\times}$ denote the set of
$(\sigma,\pmb{p})\in \hat{\mathcal{P}}$ such that $p_1,\ldots,p_r$ has
at least one intersection, and
let $\hat{\mathcal{N}}$ denote the complement
of $\hat{\mathcal{P}}^\times$
in $\hat{\mathcal{P}}.$ 
For any subset $S$ of $\hat{\mathcal{P}}$
let $G(S)$ denote the corresponding
generating functions.
Then we have
$$
G(\hat{\mathcal{P}})=G(\hat{\mathcal{N}})
+G(\hat{\mathcal{P}}^\times).
$$

We will show that $G(\hat{\mathcal{P}}^\times)=0.$
The argument is similar the one in \cite{St2}.
We construct a sign-reversing involution $(\sigma,\pmb{p})\mapsto
(\sigma',\pmb{p}')$ on the set $\hat{\mathcal{P}}^{\times}$ 
such that $w(\pmb{p})=w(\pmb{p}').$
The existence of such involution 
implies $G(\hat{\mathcal{P}}^\times)=0.$

Let $v$ be a vertex on a path $p_i.$
We denote by $p_i(\to v)$ and $p_i(v\to)$
the subpaths of $p_i$ from $u_i$ to $v$ and $v$ to the end point. 
We say a vertex $v$ is an {\it intersection point\/}
of the $r$-tuple $\pmb{p}=(p_1,\ldots,p_r)$
if there are at least one pair of paths
that intersect at $v.$ 
Suppose there is an intersecting point
of $\pmb{p}=(p_1,\ldots,p_r).$
There is an index $i$ that satisfies the following:
\begin{itemize}
\item $p_i$ intersect $p_{i+1}$ at a vertex $v,$
\item there are no intersection points of $\pmb{p}$
on $p_i(\to v)$ and $p_{i+1}(v\to).$ 
\end{itemize}
The existence of such an index $i$ can be shown by 
introducing a total order of the vertices (cf. \cite{St2}).

Define paths 
$p_i'=p_i(\to v)p_{i+1}(v\to ),\,p_{i+1}'=p_{i+1}(\to v)p_i(v\to ),$
and set $p_k'=p_k$ for $k\ne i,i+1.$
Set $\pmb{p}'=(p_1',\ldots,p'_r).$ 
Also denote by $\sigma'$ the perfect matching
obtained by interchanging $u_i$ and $u_{i+1}$ in $\sigma.$ 
Then we claim that
the $r$-tuple $\pmb{p}'$
is $\sigma'$-admissible.
If $(u_k,u_l)$ in $\sigma'$, then
the paths $p'_k$ and $p'_l$ do not intersect.
For the proof this fact we refer to the proof of Theorem 3.1 in \cite{St2}.
Next we have to show that 
each path $p_i'\;(1\leq i\leq r)$ ends at the right region
specified by the perfect matching $\sigma'.$
Suppose that $p_i$ ends at a point in $J_1$
and $p_{i+1}$ ends at a point in $J_2.$
There are indices $k,l$ such that $(u_i,u_{k}),\,(u_l,u_{i+1})\in \sigma$
with $i<k,\,l< i+1.$
Note that we have $l<i$ and $i+1<k$
since $\{i,k\}$ and $\{l,i+1\}$ are disjoint. 
Now by construction, the path $p'_i$ ends at a point in $J_2$
and $p'_{i+1}$ ends at a point in $J_1.$
Since the paths $p'_k\;(k\ne i,i+1)$ are not changed
we see that all the paths $p_1',\ldots,p'_r$ end at the right regions.
Other cases are left for the reader.
Thus we have an involution $(\sigma,\pmb{p})\mapsto (\sigma',\pmb{p}')$
on $\mathcal{P}^{\times}$ with the desired property. 

Let $\hat{\mathcal{N}}_0$ denote the 
subset of $\hat{\mathcal{N}}$ consisting of the pairs
$(\sigma_0,\pmb{p})$ such that $\pmb{p}\in \mathcal{N}_\mu^{I\!I}(\lambda).$
Obviously we have $G(\hat{\mathcal{N}}_0)=\mathcal{W}_\mu^{I\!I}(\lambda).$
It remains to show that 
we can delete all of the terms 
corresponding to the complementary 
set $\hat{\mathcal{N}}\setminus \hat{\mathcal{N}}_0.$
To show this we define a sign-reversing involution 
on the set $\hat{\mathcal{N}}\setminus \hat{\mathcal{N}}_0$
by 
$(\sigma,\pmb{p})\mapsto (\sigma^{*},\pmb{p}),$
where $\sigma^{*}$ is defined in Lemma \ref{invo}.
Note that the condition $I_\sigma=I_{\sigma^{*}}$ 
insure that the $r$-tuple $\pmb{p}$ is $\sigma^{*}$-admissible.
\hfill $\square$

\newpage 

\subsection{Bijection between 
$\mathcal{N}_\mu^{*}(\lambda)$ and $\mathcal{E}_\mu^{*}(\lambda)$ }
We will establish a bijection between 
$\mathcal{N}_\mu^{*}(\lambda)$ and $\mathcal{E}_\mu^{*}(\lambda).$
The following is an example of a tuple of
non-intersecting paths with $\lambda=(4,3,2,1)$
and $\mu=\rho_4$ and the corresponding shifted EYD:
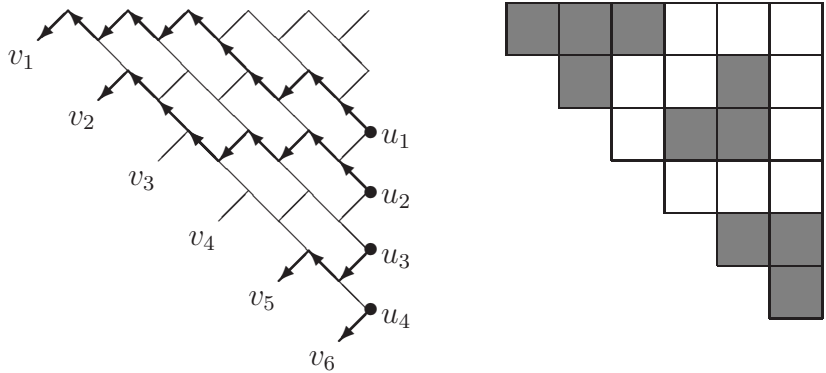
\begin{center}
  \begin{picture}(90,65)
  \thicklines
  \thinlines
  \put(25,60){\line(1,-1){50}}
  \put(35,60){\line(1,-1){40}}
  \put(45,60){\line(1,-1){30}}
  \put(55,60){\line(1,-1){20}}
  \put(65,60){\line(1,-1){10}}
  \multiput(20,55)(10,0){6}{\line(1,1){5}}
  \multiput(30,45)(10,0){5}{\line(1,1){5}}
  \multiput(40,35)(10,0){4}{\line(1,1){5}}
  \multiput(50,25)(10,0){3}{\line(1,1){5}}
  \multiput(60,15)(10,0){2}{\line(1,1){5}}
  \multiput(70,5)(10,0){1}{\line(1,1){5}}
  \put(15,51){$v_1$}
  \put(25,41){$v_2$}
  \put(35,31){$v_3$}
  \put(45,21){$v_4$}
  \put(55,11){$v_5$}
  \put(65,1){$v_6$}
  \put(77,38){$u_1$}
  \put(77,28){$u_2$}
  \put(77,18){$u_3$}
  \put(77,8){$u_4$}
  \put(74,9){$\bullet$}
  \put(74,19){$\bullet$}
  \put(74,28.5){$\bullet$}
  \put(74,38.5){$\bullet$}
  \thicklines
  \multiput(75,40)(-5,5){2}{\vector(-1,1){5}}
  \multiput(75,30)(-5,5){2}{\vector(-1,1){5}}
  \multiput(25,60)(10,0){3}{\vector(-1,-1){5}}
  \multiput(30,55)(10,0){3}{\vector(-1,1){5}}
  \multiput(60,45)(-5,5){2}{\vector(-1,1){5}}
  \put(65,50){\vector(-1,-1){5}}
  \multiput(55,40)(10,0){2}{\vector(-1,-1){5}}
  \multiput(50,35)(10,0){2}{\vector(-1,1){5}}
  \multiput(50,35)(-5,5){3}{\vector(-1,1){5}} 
  \put(35,50){\vector(-1,-1){5}}
  \multiput(65,20)(10,0){2}{\vector(-1,-1){5}}
  \multiput(70,15)(10,0){1}{\vector(-1,1){5}}
  \put(75,10){\vector(-1,-1){5}}
  \end{picture}
\setlength{\unitlength}{0.7mm}
  \begin{picture}(90,65)
  \put(5,60){\Cell}
  \put(15,60){\Cell}
  \put(25,60){\Cell}
  \put(15,50){\Cell}
  \put(45,50){\Cell}
  \put(35,40){\Cell}
  \put(45,40){\Cell}
  \put(55,20){\Cell}
  \put(45,20){\Cell}
  \put(55,10){\Cell}
  \put(5,70){\line(1,0){60}}
  \put(5,60){\line(1,0){60}}
  \put(15,50){\line(1,0){50}}
  \put(25,40){\line(1,0){40}}
  \put(35,30){\line(1,0){30}}
  \put(45,20){\line(1,0){20}}
  \put(55,10){\line(1,0){10}}
  \put(5,60){\line(0,1){10}}
  \put(15,50){\line(0,1){20}}
  \put(25,40){\line(0,1){30}}
  \put(35,30){\line(0,1){40}}
  \put(45,20){\line(0,1){50}}
  \put(55,10){\line(0,1){60}}
  \put(65,10){\line(0,1){60}}
  \end{picture}
\end{center}

Let $\pmb{p}=(p_1,\ldots,p_r)\in\mathcal{N}_\mu^{*}(\lambda).$ 
We define some subsets of $D_\mu'$ by 
$C(\pmb{p})=\bigcup_{i=1}^{r} C(p_{i}),$
where 
$
C(p_{i})=\{v\in D'_\mu\;|\;p_{i}\;
\mbox{goes through the edge from}\;
v\;\mbox{to}\;v+\pmb{a}\}.
$
First we claim the following.
\begin{prop}
Let $\pmb{p}\in \mathcal{N}_\mu^{*}(\lambda).$ Then
$C(\pmb{p})$ is an element of $\mathcal{E}_\mu^{*}(\lambda).$
\end{prop}
{\it Proof.} We use
induction on the energy of $C(\pmb{p}).$
Suppose $C(\pmb{p})$ has zero energy.@
Then we have $C(\pmb{p})=D_{\lambda}'$ and 
the Proposition is true.
Let $\pmb{p}^{0}=(p_{1}^{0},\ldots,p_{r}^{0})$ denote the 
corresponding $r$-tuple with $D_{\lambda}'$, the ground state.
If $C$ has an energy $>0$, then there is a path $p_{i}$ such 
that $p_{i}\ne p_{i}^{0}.$ We take the smallest index $i.$
Let $v$ be the `last' vertex in $C(p_{i})$  
that differs from those in $C(p_{i}^{0}).$
First consider the case that $v$ is not on the main diagonal,
i.e., $v=(i,j)$ with $i<j.$ 
Put $x_{1}=v-(1,1),\;x_{2}=v-(1,0),\;x_{3}=(0,1).$
Then we see that there are no paths 
that go through $x_{1}$ and $x_{2}$
from the way of choosing the vertex $v.$
Next we claim that $x_{3}$ does not belong to $C(\pmb{p}).$
This is clear because the path $p_{i}$ goes straight through 
the vertex $x_{3}$ to NW direction.
Now we deform the path $p_{i}$ 
to get $p_{i}'$ as in the figure,
and set $\pmb{p}'=(p_{1}',\ldots,p_{r}')$ 
with $p_{j}'=p_{j}$ for $j\ \ne i.$
Then we have $\pmb{p}'\in \mathcal{N}_{\mu}^{*}(\lambda)$
and $C(\pmb{p})$ is obtained from an elementary excitation
from $C(\pmb{p}').$
By the hypothesis of induction we have $C(\pmb{p}')\in \mathcal{E}_{\mu}^{*}(\lambda).$
So we have $C(\pmb{p})\in \mathcal{E}_{\mu}^{*}(\lambda).$
\setlength{\unitlength}{0.8mm}
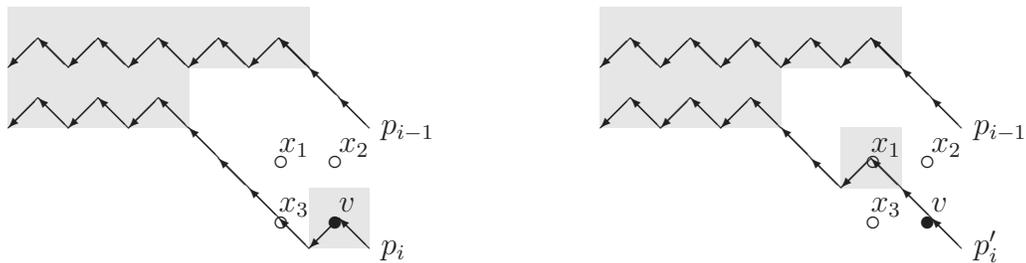
\begin{center}
  \begin{picture}(70,45)
  \put(5,35){$\Cellb$}
  \put(15,35){$\Cellb$}
  \put(25,35){$\Cellb$}
  \put(35,35){$\Cellb$}
  \put(45,35){$\Cellb$}
  \put(5,25){$\Cellb$}
  \put(15,25){$\Cellb$}
  \put(25,25){$\Cellb$}
  \put(55,5){$\Cellb$}
  \put(67,24){$p_{i-1}$}
  \multiput(10,30)(10,0){3}{\vector(-1,-1){5}}
  \multiput(15,25)(10,0){3}{\vector(-1,1){5}}
  \multiput(55,5)(-5,5){4}{\vector(-1,1){5}}
  \put(60,10){\vector(-1,-1){5}}
  \put(65,5){\vector(-1,1){5}}
  \put(67,4){$p_{i}$}
  \multiput(65,25)(-5,5){3}{\vector(-1,1){5}}
  \multiput(10,40)(10,0){5}{\vector(-1,-1){5}}
  \multiput(15,35)(10,0){5}{\vector(-1,1){5}}
  \put(49,8){$\circ$} 
  \put(49,18){$\circ$} 
  \put(58,18){$\circ$} 
  \put(58,8){$\bullet$} 
  \put(50,11){$x_{3}$} 
  \put(50,21){$x_{1}$} 
  \put(60,21){$x_{2}$} 
  \put(60,11){$v$} 
\end{picture}
\hspace{2cm}
  \begin{picture}(70,45)
  \put(5,35){$\Cellb$}
  \put(15,35){$\Cellb$}
  \put(25,35){$\Cellb$}
  \put(35,35){$\Cellb$}
  \put(45,35){$\Cellb$}
  \put(5,25){$\Cellb$}
  \put(15,25){$\Cellb$}
  \put(25,25){$\Cellb$}
  \put(45,15){$\Cellb$}
  \put(67,24){$p_{i-1}$}
  \multiput(10,30)(10,0){3}{\vector(-1,-1){5}}
  \multiput(15,25)(10,0){3}{\vector(-1,1){5}}
  \multiput(45,15)(-5,5){2}{\vector(-1,1){5}}
  \put(50,20){\vector(-1,-1){5}}
  \multiput(65,5)(-5,5){3}{\vector(-1,1){5}}
  \put(67,4){$p_{i}'$}
  \multiput(65,25)(-5,5){3}{\vector(-1,1){5}}
  \multiput(10,40)(10,0){5}{\vector(-1,-1){5}}
  \multiput(15,35)(10,0){5}{\vector(-1,1){5}}
  \put(49,8){$\circ$} 
  \put(49,18){$\circ$} 
  \put(58,18){$\circ$} 
  \put(58,8){$\bullet$} 
  \put(50,11){$x_{3}$} 
  \put(50,21){$x_{1}$} 
  \put(60,21){$x_{2}$} 
  \put(60,11){$v$} 
\end{picture}
\end{center}

Next we consider the case that $v$ is on the main diagonal, i.e., $v=(i,i).$
It is enough to consider type $I\!I$ case.
Put $x_{1}=v-(1,1),\;x_{2}=v-(2,2),\;x_{3}=v-(1,0),\,x_{4}=v-(2,1).$
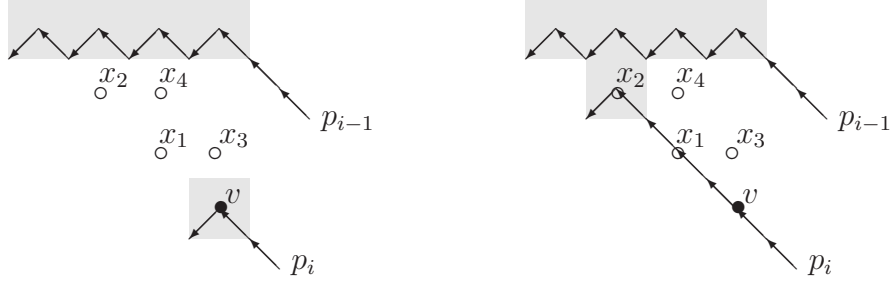
\setlength{\unitlength}{0.8mm}
\begin{center}
  \begin{picture}(70,60)
  \put(5,45){$\Cellb$}
  \put(15,45){$\Cellb$}
  \put(25,45){$\Cellb$}
  \put(35,45){$\Cellb$}
  \put(35,15){$\Cellb$}
  \put(57,34){$p_{i-1}$}
  \multiput(50,10)(-5,5){2}{\vector(-1,1){5}}
  \put(40,20){\vector(-1,-1){5}}
  \put(52,10){$p_{i}$}
  \multiput(55,35)(-5,5){2}{\vector(-1,1){5}}
  \multiput(10,50)(10,0){4}{\vector(-1,-1){5}}
  \multiput(15,45)(10,0){4}{\vector(-1,1){5}}
  \put(38,28){$\circ$} 
  \put(29,28){$\circ$} 
  \put(29,38){$\circ$} 
  \put(19,38){$\circ$} 
  \put(39,19){$\bullet$} 
  \put(40,31){$x_{3}$} 
  \put(30,31){$x_{1}$} 
  \put(30,41){$x_{4}$} 
  \put(20,41){$x_{2}$} 
  \put(41,21){$v$} 
\end{picture}
\hspace{1cm}
  \begin{picture}(70,60)
  \put(5,45){$\Cellb$}
  \put(15,45){$\Cellb$}
  \put(25,45){$\Cellb$}
  \put(35,45){$\Cellb$}
  \put(15,35){$\Cellb$}
  \put(57,34){$p_{i-1}$}
  \multiput(50,10)(-5,5){6}{\vector(-1,1){5}}
  \put(20,40){\vector(-1,-1){5}}
  \put(52,10){$p_{i}$}
  \multiput(55,35)(-5,5){2}{\vector(-1,1){5}}
  \multiput(10,50)(10,0){4}{\vector(-1,-1){5}}
  \multiput(15,45)(10,0){4}{\vector(-1,1){5}}
  \put(38,28){$\circ$} 
  \put(29,28){$\circ$} 
  \put(29,38){$\circ$} 
  \put(19,38){$\circ$} 
  \put(39,19){$\bullet$} 
  \put(40,31){$x_{3}$} 
  \put(30,31){$x_{1}$} 
  \put(30,41){$x_{4}$} 
  \put(20,41){$x_{2}$} 
  \put(41,21){$v$} 
\end{picture}
\end{center}
There are no paths going though $x_{1}$ and $x_{2}$
since we have the restriction
of the end point of the paths.
Next we claim that $x_{3},x_{4}$ do not belong to $C(\pmb{p}).$
This is clear if $i=1.$
Suppose the upper path $p_{i-1}$ exists i.e. $i>1$ 
and $x_{j}\in C(p_{i-1})$ for $i=3$ or $4$, then  
the path $p_{i-1}$ must go through $x_{1}$ or $x_{2}.$
This contradict the above assertion. 
Thus in particular $x_{i}$ for $1\leq i\leq 4$ are not in $C(\pmb{p}).$
Now we deform the path $p_{i}$ to $p_{i}'$ 
as indicated in figure
and set $\pmb{p}'=(p_{1}',\ldots,p_{r}')$ 
with $p_{j}'=p_{j}$ for $j\ \ne i.$  
Then by the same argument of the preceding case,
we have $C(\pmb{p})\in \mathcal{E}_{\mu}^{I\!I}(\lambda).$
\hfill $\square$

\bigskip

We shall describe the inverse map of $\pmb{p}\mapsto C(\pmb{p}).$
Let $C$ be an element of $\mathcal{E}_\mu^{*}(\lambda).$
For each $(i,j)\in C,$ we define its {\it layer number} 
as follows.
Let $L_{m}$ denote the diagonal line $\{(i,j)\in \Z^{2}\,|\,j-i=m\}.$ 
If $(i,j)\in L_{m}$ and 
$$
C\cap L_{m}=\{(i_{1},j_{1}),\ldots,(i_{s},j_{s})\}
$$
with $i_{1}<\cdots<i_{s},$ and $(i,j)=(i_{k},j_{k}),$
then we define the layer number of $(i,j)$ to be $k.$ 
Let $b_1,\ldots,b_m\in D'_\mu$ be the boxes in $C$
with the layer number $k.$
It is obvious that there exists
a unique path starting from $u_k$
and going through all the points $b_1,\ldots,b_m.$
Let us denote the path by $p_k.$  
Then form an $r$-tuple $\pmb{p}_C=(p_1,\ldots,p_r).$
\begin{lem} 
$\pmb{p}_C$ is non-intersecting.
\end{lem}
{\it Proof.} Clear from the construction.
\hfill$\square$

\bigskip

Obviously if $C\in \mathcal{E}_\mu^{I\!I}(\lambda),$
then $\pmb{p}_C\in \mathcal{N}_\mu^{I\!I}(\lambda).$
For any subset $C$ of $D_\mu'$, we define 
its weight by $w(C)=\prod_{(i,j)\in C}a_{ij}.$

\begin{prop}\label{bijEN} Let $\lambda\leq \mu$ be strict partitions.
There exists a weight preserving
bijection between $\mathcal{N}_\mu^I(\lambda)$ 
and $\mathcal{E}_\mu^I(\lambda).$
Moreover, this map induces a weight preserving
bijection between $\mathcal{N}_\mu^{I\!I}(\lambda)$ 
and $\mathcal{E}_\mu^{I\!I}(\lambda).$
Thus we have
$$
\sum_{C\in \mathcal{E}_\mu^{*}(\lambda)}w(C)
=\sum_{\pmb{p}\in \mathcal{N}_\mu^{*}(\lambda)}w(\pmb{p}),
$$
for $*=I,I\!I.$ 
\end{prop}

Let $w\leq v\in W^{P},$ 
and $\lambda\leq \mu\in \mathcal{SP}_{\!\rho_n}$  the
corresponding strict partitions. 
Then our combinatorial formula (Theorem \ref{BCD}) 
reads 
$$
[X_\lambda]|_{v} =\sum_{C\in\mathcal{E}^{*}_\mu(\lambda)}w(C),
$$
where we chose the weights $a_{ij}=\beta_{i,j}$ which are given in Lemma \ref{betaBCD}.
Now from Propositions \ref{LPMPf} and \ref{bijEN}, we have
$$
[X_\lambda]|_{\mu}
=\mathrm{Pf}\left([X_{\lambda_i,\lambda_j}]|_\mu
\right),
$$
which is equal to 
$\left(\mathrm{Pf}\,[X_{\lambda_i,\lambda_j}]\right)|_\mu$
since the map $\alpha\mapsto \alpha|_{v}$ is a ring homomorphism.
By virtue of the injectivity of the localization map,
we have the equivariant Giambelli formula.

\begin{small}
{\scshape Department of Applied Mathematics,
Okayama University of Science,
Okayama 700-0005, JAPAN}
\end{small}

{\textit{email address}: \tt{ike@xmath.ous.ac.jp}}

\bigskip

\begin{small}
{\scshape Faculty of Education,
Okayama University,
Okayama 700-8530, JAPAN}
\end{small}

\textit{email address}: \tt{rdcv1654@cc.okayama-u.ac.jp}

\end{document}